\documentclass[a4paper,11pt,reqno,twoside]{amsart}

\usepackage{amsmath}
\usepackage{amsfonts}
\usepackage{amssymb}
\usepackage{amsthm}
\usepackage{color}
\usepackage{ifpdf}
\usepackage{array}
\usepackage{url}
\usepackage{multirow,bigdelim}
\usepackage{mathrsfs}
\usepackage[all]{xy}
\usepackage{graphicx}

\newtheorem{theorem}{Theorem}[section]
\newtheorem{proposition}[theorem]{Proposition}
\newtheorem{lemma}[theorem]{Lemma}
\newtheorem{corollary}[theorem]{Corollary}

\numberwithin{equation}{section}

\newcommand*{\C}{\mathbb{C}}
\newcommand*{\R}{\mathbb{R}}
\newcommand*{\D}{\mathscr{D}}
\newcommand{\Ha}{\mathcal H(\mathsf A_a)}

\newcommand{\comment}[1]{}
\title[The Local Weil Explicit Formula]%
      {The Local Weil Explicit Formula \\ and Related Structures}
\author[M. Suzuki]{Masatoshi Suzuki}
\date{Version of \today}
\subjclass[2020]{
11M26, 
42A82, 
46E22, 
47B25
}
\keywords{
Weil quadratic form;
local Weil explicit formula;
canonical systems;
de Branges spaces;
screw functions%
}
\AtBeginDocument{%
\begin{abstract}
We study the reproducing kernel Hilbert spaces constructed from the local
Weil quadratic form associated with the Riemann zeta function.
Wherever this form is positive definite, we prove a local explicit formula
expressing the Weil distribution in terms of the discrete spectra of
self-adjoint extensions of the differentiation operator.
We characterize the Riemann hypothesis by the uniform boundedness
of the diagonal values of the reproducing kernels.
We then study the chain of de Branges spaces and recover the Hamiltonian
of the corresponding canonical system explicitly from these kernels
under local positivity, without assuming the Riemann hypothesis.
Finally, assuming the Riemann hypothesis, we study the asymptotic
behavior of solutions of this system.
\end{abstract}
\maketitle
}
\begin{document}

\section{Introduction}

The equality 
\begin{equation} \label{eq_101}
\sum_{p^k \leq x} \log p
= x 
- \sum_{\rho} \frac{x^{\rho}}{\rho}
- \log 2\pi - \frac{1}{2}\log(1-x^{-2})
\quad (x>1,\ x\ne p^k)
\end{equation}
relating the primes to the nontrivial zeros of the Riemann zeta function $\zeta(s)$
plays an essential role in the proof of the prime number theorem
using the zeta function.
Here $\sum_\rho$ runs over the nontrivial zeros of $\zeta(s)$,
counted with multiplicity, and denotes the limit of the sums
over $|\Im(\rho)|\le T$ as $T\to\infty$ (see \cite[(12.1)]{MV07}).
In analytic number theory, equalities of this kind relating prime powers
to zeros of zeta functions are known as explicit formulas.
The following \emph{Weil explicit formula} generalizes this equality by
incorporating test functions satisfying suitable conditions.
For $u \in C_c^\infty(\R)$, define
\begin{equation} \label{eq_102}
\aligned 
W(u)
& := \int_{-\infty}^{\infty} u(x) (e^{x/2}+e^{-x/2}) dx 
 - \sum_{n=1}^{\infty} \frac{\Lambda(n)}{\sqrt{n}} u(\log n) 
- \sum_{n=1}^{\infty} \frac{\Lambda(n)}{\sqrt{n}} u(-\log n)  \\
& \quad 
- (\log 4\pi + C_0) u(0)
- \int_{0}^{\infty} \left\{ u(x) + u(-x)-2e^{-x/2} u(0)\right\} \frac{e^{x/2}dx}{e^{x}-e^{-x}}, 
\endaligned 
\end{equation}
where $\Lambda$ is the von Mangoldt function,
\[
\Lambda(n) =
\begin{cases}
~\log p & \text{if, for some prime $p$, $n=p^k$ ($k\ge1$)}, \\
~0 & \text{otherwise}, 
\end{cases}
\]
and $C_0$ is the Euler--Mascheroni constant.
The linear functional $u \mapsto W(u)$ 
is called the Weil distribution 
(see \cite[proof of Proposition~2.1]{Su23JNT}
and \cite[\S1.1]{Su26a}).

Let $\xi(s):=(1/2)s(s-1)\pi^{-s/2}\Gamma(s/2)\zeta(s)$ be the Riemann xi-function. 
For the functional equation and basic properties of the zeros 
of $\xi(s)$, see \cite[\S10.1, \S10.2]{MV07}. 
Let $\Gamma$ denote the set of zeros of the function
$z \mapsto \xi(1/2-iz)$, and, for each $\gamma \in \Gamma$,
let $m_\gamma$ be the multiplicity of the zero $\gamma$.
The elements of $\Gamma$ need not be real.
We use the Fourier transform
\[
\widehat{u}(z)
:= \int_{-\infty}^{\infty}u(x)e^{izx}\,dx.
\]
Then, for every $u\in C_c^\infty(\R)$,
we have
\begin{equation} \label{eq_103}
\sum_{\gamma\in\Gamma}m_\gamma\,\widehat{u}(\gamma)
= W(u).
\end{equation}
This is the Weil explicit formula
(see \cite[Theorem~12.13]{MV07}
and \cite[(5.15)]{Su23}). 
The class of admissible test functions $u$ can be extended to include
discontinuous functions, such as indicator functions of intervals,
and the classical formula \eqref{eq_101} can be regarded as
a special case. However, we first work with $u \in C_c^\infty(\R)$.
\medskip

The Weil explicit formula has many applications in number theory, 
including a proof of the prime number theorem and a criterion for the Riemann hypothesis.
The Riemann hypothesis (RH), which asserts that every nontrivial zero
of the Riemann zeta function lies on the critical line $\Re(s)=1/2$,
is equivalent to the nonnegativity of the Weil quadratic form
defined by the Weil distribution,
\begin{equation*}
Q_W(u,f):=W(u \ast \widetilde{f}), \qquad Q_W(u):=Q_W(u,u),
\end{equation*}
for every $u \in C_c^\infty(\R)$. 
This is known as Weil's positivity criterion
\cite{We52,Yo92}. 
Here
\[
(u\ast f)(x) := \int_{-\infty}^{\infty}
u(y)f(x-y)\, dy, \qquad 
\widetilde{u}(x) := \overline{u(-x)}. 
\]

At the time of writing, RH was still an open problem,
and hence it was unknown whether $Q_W$ is nonnegative.\footnote{
Shortly after the completion of this work,
L.~Alp\"{o}ge announced a proof of RH in Lean
on November 20, 2026.
This development does not affect the results of the present paper.
}
However, it is known that $Q_W(u)>0$ if the support of the test function
$u \not=0$ is restricted to a sufficiently small interval.
Indeed, Yoshida~\cite[proof of Lemma 2]{Yo92} showed that, for all sufficiently
small $a>0$, there exists $c=c_a>0$ such that
$Q_W(u)\geq c \|u\|_{L^2}^2$ for every $u\in C_c^\infty(-a,a)$.
Here $C_c^\infty(-a,a)$ denotes the space of $C^\infty$ functions
whose support is contained in $(-a,a)$.
The restriction of $Q_W$ to $C_c^\infty(-a,a)$ is closable
as a lower bounded quadratic form on $L^2(-a,a)$ \cite[\S3]{Su26a}.
We denote its closure by $Q_W^a$ and the unique self-adjoint
operator associated with $Q_W^a$ by $\mathsf A_a$. 
We write $Q_W^a>0$ if there exists $c>0$ such that
$Q_W^a(u,u)\ge c\|u\|_{L^2}^2$ throughout the form domain.
The preceding estimate extends to $Q_W^a$.

The purpose of this paper is to obtain a local version of the Weil explicit
formula from this local positivity and to describe the structure of the
associated family of reproducing kernel Hilbert spaces.
The results develop the study initiated in~\cite{Su26a}.
\medskip

We now describe the local explicit formula.
Assume first that $Q_W^a$ is positive, and let $\Ha$ be the Hilbert space
obtained by completing $C_c^\infty(-a,a)$ with respect to $Q_W^a$.
The space of Fourier transforms of its elements,
\[
\mathcal B_a:=\widehat{\Ha},
\]
is a reproducing kernel Hilbert space; we denote its reproducing kernel
by $J_a(w,z)$.
The differentiation operator $\D_a=i\,d/dx$ on $\Ha$, with domain
$C_c^\infty(-a,a)$, has a family of self-adjoint extensions
$\overline{\mathscr D}_{a,\theta}$ parametrized by $\theta \in [0,2\pi)$.
Each extension has compact resolvent (Proposition~\ref{prop_202}).
Let $\Gamma_{a,\theta}$ denote its set of eigenvalues.
We then have the following result.

\begin{theorem}[Local Weil explicit formula] \label{thm_1}
Fix $a>0$ such that $Q_W^a>0$.
Then, for every $\theta\in[0,2\pi)$ and
$u \in C_c^\infty(-a,a)$, we have 
\begin{equation}\label{eq_104}
 W(u)=
 \sum_{\gamma\in\Gamma_{a,\theta}}
       \frac{\widehat{u}(\gamma)}{J_a(\gamma,\gamma)},
\end{equation}
and the series on the right converges absolutely.
\end{theorem}

For a fixed $a>0$ such that $Q_W^a>0$, the left-hand side of
\eqref{eq_104} is independent of the phase $\theta$ specifying
the self-adjoint extension of $\D_a$, whereas the right-hand side
is a sum over the set of eigenvalues $\Gamma_{a,\theta}$ corresponding
to the phase $\theta$.
Moreover, for $u\in C_c^\infty(-a,a)$,
only finitely many terms in the sums in \eqref{eq_102} are nonzero,
and the integrals are over bounded intervals
apart from the term involving $u(0)$.
Thus Theorem~\ref{thm_1} represents the Weil distribution,
determined by arithmetic data from a bounded range,
as a sum over a discrete spectrum.

These spectral representations are not merely obtained separately
for each $a>0$; they are compatible as $a$ varies.
Fix $R>0$ and suppose that $Q_W^a>0$ for every $0<a<R$.
If $0<a<b<R$, then
$C_c^\infty(-a,a)\subset C_c^\infty(-b,b)$,
and the restrictions of the Weil distribution to these spaces agree
on the smaller space.
The inclusions therefore induce isometric embeddings
$
\mathcal H(\mathsf A_a)\hookrightarrow\mathcal H(\mathsf A_b)
$ and
$
\mathcal B_a\hookrightarrow\mathcal B_b
$.
On the other hand, for a fixed $u\in C_c^\infty(-a,a)$, 
the local Weil explicit formula holds for every $a\le b<R$.
Although the set of eigenvalues $\Gamma_{b,\theta}$ and the weights
$J_b(\gamma,\gamma)^{-1}$ vary with $b$, the resulting spectral
representations always equal $W(u)$.

Thus the Hilbert spaces obtained from the local Weil quadratic form
constitute a chain under isometric embeddings, and the spectral
representations at each point of the chain are related 
through the same Weil distribution.
This structure holds wherever the local Weil quadratic form is positive,
\emph{without} assuming RH.

\medskip

In fact, each reproducing kernel Hilbert space $\mathcal B_a$ is a
de Branges space, and the isometric embeddings above define a chain
of de Branges spaces.
By the general theory of de Branges spaces
\cite[Theorem 40 and Problem 153]{dB}, such a chain can be described
by a system of first-order linear ordinary differential equations,
parametrized by $z \in \C$, called a \emph{canonical system}.
In the present setting, this system takes the following form in terms
of the real entire functions $A_a(z), B_a(z)$ obtained from a generator
of $\mathcal B_a$ with the normalization specified in \eqref{eq_502}:
\begin{equation}\label{eq_105}
\begin{gathered}
J\frac{\partial}{\partial a}
\begin{pmatrix}A_a(z)\\B_a(z)\end{pmatrix}
=
-z\begin{pmatrix}\gamma(a)&0\\0&\gamma(a)^{-1}\end{pmatrix}\begin{pmatrix}A_a(z)\\B_a(z)\end{pmatrix}, 
\qquad
\begin{pmatrix}A_0(z)\\B_0(z)\end{pmatrix}
=
\begin{pmatrix}1\\0\end{pmatrix},\\[6pt]
J=\begin{pmatrix}0&-1\\1&0\end{pmatrix},\quad z\in\C,\quad\text{a.e. }a\in(0,R).
\end{gathered}
\end{equation}
The existence of a canonical system is guaranteed by the general theory, 
but the coefficient $\gamma(a)$ is not explicitly determined in general.
\medskip

We show that the second row of the fundamental matrix of a canonical system
associated with \eqref{eq_105} is $(-B_a(z),A_a(z))$, where the real
entire functions $A_a(z),B_a(z)$ are determined by the reproducing kernel
$J_a$ of $\mathcal B_a$ (Proposition~\ref{prop_505}).
This identifies the local spectra corresponding to the phases $0,\pi$
with the zeros of entries of the fundamental matrix of the canonical system.
The local Weil explicit formula and the theorem of Krein and Langer
on extensions of \emph{screw functions} (Theorem~\ref{thm_404})
play a central role in this identification.
This identification also enables us to recover
the coefficient $\gamma(a)$ of the canonical system
\eqref{eq_105} explicitly from the reproducing kernel $J_a$. 

\begin{theorem} \label{thm_2}
Suppose that $Q_W^a>0$ for every $0<a<R$, and let $J_a(w,z)$ be
the reproducing kernel of $\mathcal B_a$.
Then the coefficient $\gamma(a)$ in the canonical system
\eqref{eq_105} is given by
\begin{equation} \label{eq_106}
\gamma(a)=
\frac{J_a(i,0)}{J_a(i,i)}
\frac{\partial}{\partial a}J_a(i,0)
\quad\text{a.e. on $(0,R)$}.
\end{equation}
\end{theorem}

So far we have considered only those $a>0$ for which $Q_W^a>0$.
None of these results requires RH.

Under RH, we have $Q_W^a>0$ for every $a>0$ \cite[\S2.2]{Su26b},
so the chain of Hilbert spaces described above extends over the entire
positive half-line.
Moreover, by \cite[Theorem~1.1]{Su26b}, $\Gamma_{a,\pi}$ converges
to $\Gamma$ as $a\to\infty$. 
If $\gamma_a\in\Gamma_{a,\pi}$ is chosen suitably so that it converges
to $\gamma\in\Gamma$, then
$J_a(\gamma_a,\gamma_a)\to 1/m_\gamma$ $(a\to\infty)$.
We prove this in Section~\ref{section_2}.
Thus, as $a\to\infty$, the eigenvalues and their weights 
in the local Weil explicit formula recover the zeta zeros and their multiplicities in the usual Weil explicit formula.

If RH is false, on the other hand, \cite[Theorem 1.3]{Su26a} implies
that there exists $a_\ast>0$ such that $Q_W^a>0$ for $0<a<a_\ast$,
whereas at $a=a_\ast$ there is at least one nonzero function $u$
in the domain of $Q_W^{a_\ast}$ for which $Q_W^{a_\ast}(u)=0$.
Examining this difference in the behavior of the chain $\mathcal B_a$
through its reproducing kernels leads to the following result.

\begin{theorem}\label{thm_3}
RH holds if and only if the following uniform boundedness condition holds:
for every $R>0$ such that $Q_W^a>0$ for all $0<a<R$, there exists
$C_R>0$ such that
\[
 J_a(t,t)\le C_R\qquad(0<a<R,\ t\in\R).
\]
\end{theorem}

This result detects the failure of RH through the unboundedness of the
diagonal values $J_a(t,t)$ of the reproducing kernels within a bounded
range of $a>0$.
Since the construction of the reproducing kernel involves
${\mathsf A_a}^{-1}$, this unboundedness can also be interpreted as
a reformulation of the fact that $0$ is an eigenvalue of $\mathsf A_a$
at $a=a_\ast$.
(At $a=a_\ast$, the kernel $J_a(w,z)$ is not defined, but
$\mathsf A_a$ is.)
\medskip

Finally, assuming RH, we study the canonical system on the half-line.
The Hamiltonian recovered in Theorem~\ref{thm_2} gives a Krein string
of finite length and infinite mass. Using the convergence of the local
spectra in \cite[Theorems~1.1 and~1.2]{Su26b}, we determine the asymptotic 
behavior of the fundamental matrix and show that the reproducing kernel,
normalized by $J_a(0,0)$, converges to a kernel of rank one as
$a\to\infty$ (Proposition~\ref{prop_801} and Theorem~\ref{thm_802}).
\medskip

The paper is organized as follows.
In Section~\ref{section_2}, we prove the local Weil explicit formula
(Theorem~\ref{thm_1}).
In Section~\ref{section_3}, we prove the characterization
of RH by the uniform boundedness of the diagonal values of the reproducing
kernels (Theorem~\ref{thm_3}).
Section~\ref{section_4} reviews the required results on de Branges spaces,
canonical systems, and extensions of screw functions.
In Section~\ref{section_5}, we study the chain of de Branges spaces obtained
from the local Weil quadratic form and construct a particular canonical
system.
In Section~\ref{section_6}, we recover its Hamiltonian explicitly and use
it to prove Theorem~\ref{thm_2}.
In Section~\ref{section_7}, we study the relation between
the local spectral measures and the zeta zeros, extensions that give global
screw functions, and the length and mass of the corresponding strings.
Finally, in Section~\ref{section_8}, assuming RH, we determine the limits
of the length and mass and the asymptotic behavior of the fundamental
matrix and the reproducing kernel.

\section{Local Weil quadratic forms and the explicit formula} \label{section_2}

\subsection{Spectral representations for local Weil quadratic forms}

Throughout this paper, inner products are linear in the first variable.
The local Weil quadratic form $Q_W^a$ has domain
\begin{equation}\label{eq_201}
 H^{\log}(-a,a):=\left\{u\in L^2(-a,a) ~\left|~
 \int_{-\infty}^\infty |\widehat u(t)|^2\log(2+|t|)\,dt<\infty\right.\right\}
\end{equation}
and $C_c^\infty(-a,a)$ is a core
\cite[\S3]{Su26a}.
Set
\[
\lambda_a:=\inf\sigma(\mathsf A_a).
\]
Then $Q_W^a>0$ if and only if $\lambda_a>0$.
If $0<a<b$, extension by zero preserves both 
the values of $Q_W^a$ on test functions and their $L^2$ norms.
Taking the infimum of the corresponding Rayleigh quotients 
therefore gives $\lambda_a\ge\lambda_b$.

When $\lambda_a>0$, the Hilbert space
$\mathcal H(\mathsf A_a)$ defined in the introduction coincides, as a set, with
$H^{\log}(-a,a)=\mathfrak D(\mathsf A_a^{1/2})$,
and its inner product is given by
\[
\langle u,f\rangle_{\mathsf A_a}:=Q_W^a(u,f).
\]
Moreover, $u\in H^{\log}(-a,a)$ belongs to $\mathfrak D(\mathsf A_a)$ if and only if
there is a $g\in L^2(-a,a)$ such that
$Q_W^a(u,f)=\langle g,f\rangle_{L^2}$
for every $f\in H^{\log}(-a,a)$;
in this case, $\mathsf A_a u=g$.

\begin{lemma} \label{lem_201}
For $u,f\in H^{\log}(-a,a)$, substituting $u*\widetilde f$ into \eqref{eq_102}
gives absolutely convergent integrals and finite sums over prime powers, and the resulting value equals $Q_W^a(u,f)$.
\end{lemma}
\begin{proof}
By definition,
$u*\widetilde f$ is continuous.
Writing $\tau_xu(y)=u(y-x)$, we have
\[
2(u*\widetilde f)(0)
-(u*\widetilde f)(x)
-(u*\widetilde f)(-x)
=\langle\tau_xu-u,\tau_xf-f\rangle_{L^2}.
\]
Plancherel's identity and
$|e^{iv}-1|\le\min(2,|v|)$ give
\[
\begin{aligned}
\int_0^1\|\tau_xu-u\|_{L^2}^2\frac{dx}{x}
&=\frac1{2\pi}\int_{-\infty}^{\infty}|\widehat u(t)|^2
  \int_0^1|e^{ixt}-1|^2\frac{dx}{x}\,dt\\
&\ll\int_{-\infty}^{\infty}\log(2+|t|)|\widehat u(t)|^2\,dt
<\infty.
\end{aligned}
\]
The same holds for $f$.
Thus the Cauchy--Schwarz inequality with respect to $dx/x$ and
$1-e^{-x/2}=O(x)$ show that
the integral in \eqref{eq_102} converges absolutely near the origin.
Also, $u*\widetilde f$ is supported in $[-2a,2a]$
and has absolute value at most
$\|u\|_{L^2}\|f\|_{L^2}$ at every point.
Hence the sums over prime powers are finite,
and the remaining integrals converge absolutely.

The weighted norm associated with \eqref{eq_201} and the norm induced by the closed quadratic form $Q_W^a$, after adding a sufficiently large multiple of the squared $L^2$ norm, are both complete on $H^{\log}(-a,a)$ and continuously embed this space into $L^2(-a,a)$.
By the closed graph theorem, these two norms are equivalent. 
The preceding estimates show that the right-hand side of \eqref{eq_102}
defines a continuous sesquilinear form on $H^{\log}(-a,a)$.
Since $C_c^\infty(-a,a)$ is a core for the quadratic form $Q_W^a$,
its value agrees with $Q_W^a(u,f)$.
\end{proof}

The embedding of $H^{\log}(-a,a)$ into $L^2(-a,a)$ is compact
\cite[\S4.1]{Su26a}. Indeed, restricting the Fourier transform to a finite interval $[-T,T]$ defines
a Hilbert--Schmidt operator from $L^2(-a,a)$ to $L^2(-T,T)$,
and the squared $L^2$ norm over $|t|>T$ is bounded by $1/\log(2+T)$ times the integral in \eqref{eq_201}.
Consequently, the operator $\mathsf A_a$ associated with $Q_W^a$ has compact resolvent for every $a>0$.

Henceforth, assume that $\lambda_a>0$.
The representation of point evaluations in \cite[\S6.2, \S6.3]{Su26a} shows that
$\mathcal B_a=\widehat{\mathcal H(\mathsf A_a)}$
is a reproducing kernel Hilbert space.
Put $e_z(x)=e^{-izx}$ $(-a<x<a)$ and
\[
u_{z,a}:=\mathsf A_a^{-1}e_z
\in\mathfrak D(\mathsf A_a).
\]
Then its reproducing kernel $J_a$ is given by
\begin{equation}\label{eq_202}
 J_a(w,z)=\widehat{u_{\bar w,a}}(z)
 =\langle u_{\bar w,a},u_{\bar z,a}\rangle_{\mathsf A_a}.
\end{equation}
Since point evaluations are nonzero, $J_a(w,w)>0$ for every $w\in\C$.
Moreover, \eqref{eq_102} shows that reflection $u(x)\mapsto u(-x)$ is an isometry,
and the reproducing property gives $J_a(-w,-z)=J_a(w,z)$.
Define $\mathscr D_a f=if'$ on
$C_c^\infty(-a,a)\subset\mathcal H(\mathsf A_a)$.
Its closure has deficiency indices $(1,1)$,
and $\ker(\mathscr D_a^*-z)=\C u_{z,a}$
\cite[\S6]{Su26a}.
Let $\overline{\mathscr D}_{a,\theta}$ be the self-adjoint extension
with domain
$\mathfrak D(\overline{\mathscr D}_a)
\dotplus\C(u_{i,a}+e^{i\theta}u_{-i,a})$.
Here the deficiency vectors $u_{\pm i,a}$ have equal norms.
The eigenvalues of this extension are precisely the zeros of
\begin{equation} \label{eq_203}
 W(a,\theta;z)
 =(z-i)J_a(-i,z)+e^{i\theta}(z+i)J_a(i,z),
\end{equation}
all of which are real \cite[Theorem~1.5]{Su26a}.
The eigenspaces are one-dimensional, and the zeros of $W(a,\theta;\cdot)$ are also simple.
Indeed, for an eigenvalue $t$, write $u_{t,a}=u_0+c(u_{i,a}+e^{i\theta}u_{-i,a})$,
where $u_0\in\mathfrak D(\overline{\mathscr D}_a)$.
Apply $\mathscr D_a^*$ to this decomposition and take Fourier transforms.
By the continuity of point evaluations and the definition of the graph closure,
\[
\widehat{\overline{\mathscr D}_a u_0}(z)=z\widehat u_0(z).
\]
Thus the eigenvalue equation and \eqref{eq_203} give
\[
 cW(a,\theta;z)=(z-t)J_a(t,z)
\]
(see \cite[Lemma~3.2]{Su26b} for phase $\pi$).
Since $J_a(t,t)>0$, we have $c\ne0$ and hence $\partial_z W(a,\theta;t)\ne0$.

Define the weighted discrete measure on this spectrum by
\begin{equation} \label{eq_204}
 \mu_{a,\theta}:=
 \sum_{\gamma\in\Gamma_{a,\theta}}
 \frac{\delta_\gamma}{J_a(\gamma,\gamma)}.
\end{equation}

\begin{proposition} \label{prop_202}
Suppose that $\lambda_a>0$.
For each $\theta\in[0,2\pi)$, we have 
\begin{equation} \label{eq_205}
 \sum_{\gamma\in\Gamma_{a,\theta}}
  \frac{\widehat u(\gamma)
        \overline{\widehat f(\gamma)}}
       {J_a(\gamma,\gamma)}
 =Q_W^a(u,f)=W(u*\widetilde f)
\end{equation}
for every $u,f\in\mathcal H(\mathsf A_a)$. 
Here $W$ on the right denotes
the extension given by Lemma~\ref{lem_201}.
Moreover, $\overline{\mathscr D}_{a,\theta}$
has compact resolvent,
$\{u_{\gamma,a}/\sqrt{J_a(\gamma,\gamma)}\}
_{\gamma\in\Gamma_{a,\theta}}$
is an orthonormal basis of $\mathcal H(\mathsf A_a)$,
and $u\mapsto(\widehat u(\gamma))_\gamma$ is
an isometric isomorphism onto $L^2(\mu_{a,\theta})$.
Under RH, the same value equals
$\sum_{\gamma\in\Gamma}m_\gamma
\widehat u(\gamma)\overline{\widehat f(\gamma)}$.
\end{proposition}
\begin{proof}
For a test function $u \in C_c^\infty(-a,a)$, applying the lower bound for the quadratic form to $u,iu'$ gives
$\lambda_a(\|u\|_{L^2}^2+\|u'\|_{L^2}^2) \le Q_W^a(u,u)+Q_W^a(iu',iu')$.
Hence the domain of $\overline{\mathscr D}_a$,
equipped with the graph norm, embeds continuously into $H_0^1(-a,a)$. 
The embedding $H_0^1(-a,a)\hookrightarrow H^{\log}(-a,a)$
is compact.
Indeed, this follows by combining the $L^2$ convergence
given by Rellich's theorem with the uniform estimate
for the contribution from frequencies $|t|>M$ implied by
\[
\sup_{|t|>M}\frac{\log(2+|t|)}{1+t^2}
\longrightarrow 0\qquad(M\to\infty).
\]
Within the domain of each self-adjoint extension,
equipped with the graph norm,
$\mathfrak D(\overline{\mathscr D}_a)$ is closed and has codimension one. 
Thus the domain of each self-adjoint extension
also embeds compactly into $H^{\log}(-a,a)$.
Consequently, each extension has compact resolvent.
Moreover, by \cite[Theorem~1.5, \S6]{Su26a},
for each $\gamma\in\Gamma_{a,\theta}$,
the function $u_{\gamma,a}$ is
an eigenfunction of $\overline{\mathscr D}_{a,\theta}$
with eigenvalue $\gamma$.
The spectral theorem therefore shows that
\[
\left\{
\frac{u_{\gamma,a}}{\sqrt{J_a(\gamma,\gamma)}}:
\gamma\in\Gamma_{a,\theta}
\right\}
\]
is an orthonormal basis of $\mathcal H(\mathsf A_a)$.
By \eqref{eq_202}, the Fourier coefficients of each $u\in\mathcal H(\mathsf A_a)$
with respect to this basis are
\[
\left\langle u,
\frac{u_{\gamma,a}}{\sqrt{J_a(\gamma,\gamma)}}
\right\rangle_{\mathsf A_a}
=
\frac{\widehat u(\gamma)}{\sqrt{J_a(\gamma,\gamma)}}.
\]
Parseval's identity gives the first equality in \eqref{eq_205}.
The completeness of the eigenfunctions also shows that
$u\mapsto(\widehat u(\gamma))_\gamma$ is
an isometric isomorphism onto $L^2(\mu_{a,\theta})$.
Lemma~\ref{lem_201} completes the proof of
\eqref{eq_205}.

Under RH, Weil's explicit formula gives
$Q_W(u,f)=\sum_{\gamma\in\Gamma}m_\gamma
\widehat u(\gamma)\overline{\widehat f(\gamma)}$
for $u,f\in C_c^\infty(-a,a)$.
In particular, $Q_W(u)=\sum_{\gamma\in\Gamma}m_\gamma
|\widehat u(\gamma)|^2$, so
$u\mapsto(\widehat u(\gamma))_{\gamma\in\Gamma}$
is an isometry with respect to the norm defined by $Q_W^a$.
For any $u\in\mathcal H(\mathsf A_a)$,
choose $u_n\in C_c^\infty(-a,a)$ such that
$u_n\to u$ in the norm defined by $Q_W^a$.
The preceding isometry gives
\[
\sum_{\gamma\in\Gamma}m_\gamma
|\widehat u_n(\gamma)-\widehat u_m(\gamma)|^2
=Q_W^a(u_n-u_m)~\to~0.
\]
On the other hand, evaluation of the Fourier transform at each point
is continuous on $\mathcal H(\mathsf A_a)$,
so $\widehat u_n(\gamma)\to\widehat u(\gamma)$.
It follows that
$Q_W^a(u,u)=\sum_{\gamma\in\Gamma}m_\gamma
|\widehat u(\gamma)|^2$.
Hence, we have 
$Q_W^a(u,f)=\sum_{\gamma\in\Gamma}m_\gamma
\widehat u(\gamma)\overline{\widehat f(\gamma)}$ 
for all $u,f\in\mathcal H(\mathsf A_a)$. 
The sum on the right converges absolutely by the Cauchy--Schwarz inequality.
\end{proof}

\subsection{Proof of Theorem~\ref{thm_1}} 

By Proposition~\ref{prop_202},
\begin{equation} \label{eq_206}
 W(u*\widetilde f)=
 \int_{-\infty}^{\infty}\widehat u(t)\overline{\widehat f(t)}
 \,d\mu_{a,\theta}(t)
 \qquad(u,f\in C_c^\infty(-a,a)).
\end{equation}
We now extend \eqref{eq_206} to obtain the local Weil explicit formula
for every $u\in C_c^\infty(-2a,2a)$. 
To this end, we first establish a growth estimate for $\mu_{a,\theta}$. 

Choose a nonnegative even function $\psi\in C_c^\infty(-1/4,1/4)$
with $\int\psi=1$,
and, for $T\ge\max(1,a^{-1})$,
put $u_T(x)=T^{1/2}\psi(Tx)$.
For $|t|\le T$,
$|\widehat u_T(t)|^2\ge T^{-1}\cos^2(1/4)$.
On the other hand, the equivalence of the form norm and the $H^{\log}$ norm gives
\[
 Q_W^a(u_T,u_T)\ll_a
 \int_{-\infty}^{\infty}\log(2+|t|)|\widehat u_T(t)|^2\,dt
 \ll\log(2+T).
\]
Thus \eqref{eq_206} gives
$\mu_{a,\theta}([-T,T])\ll_a T\log(2+T)$.
By monotonicity,
$\mu_{a,\theta}([-T,T])\ll_a(1+T)\log(2+T)$
for all $T\ge0$. This estimate is uniform in $\theta$.

For every $u\in C_c^\infty(-2a,2a)$,
the rapid decay of $\widehat u$ therefore gives
$\int_{-\infty}^{\infty}|\widehat u|\,d\mu_{a,\theta}<\infty$,
and
\[
 u\longmapsto
 \int_{-\infty}^{\infty}\widehat u(t)\,d\mu_{a,\theta}(t)
\]
defines a continuous linear functional on $C_c^\infty(-2a,2a)$.
By \eqref{eq_102}, $W$ is also continuous on this space.
To prove the theorem, it therefore suffices to show that
$
 \operatorname{span}
 \{u*\widetilde f:u,f\in C_c^\infty(-a,a)\}
$
is dense in $C_c^\infty(-2a,2a)$
with respect to its usual topology.

Take any $u\in C_c^\infty(-2a,2a)$.
A partition of unity expresses $u$ as
a finite sum of smooth functions, each supported
in an interval $(c-\delta,c+\delta)$
with $|c|/2+\delta<a$.
It therefore suffices to consider the case in which
$u$ is supported in such an interval.

Set $s=c/2$ and
$u_s(x)=u(x+s)$.
For a real, even approximate identity $\eta_\varepsilon$,
put $f_\varepsilon(x)=\eta_\varepsilon(x+s)$.
For sufficiently small $\varepsilon>0$,
we have $u_s,f_\varepsilon\in C_c^\infty(-a,a)$ and
$ u_s\ast\widetilde f_\varepsilon
 =u\ast\eta_\varepsilon\to u$
uniformly with all derivatives, with supports contained
in a common compact set.
This proves the required density.
The two continuous linear functionals that agree in \eqref{eq_206}
thus agree on all of
$C_c^\infty(-2a,2a)$.
It follows that \eqref{eq_104} holds,
with an absolutely convergent right-hand side.
\hfill $\Box$
\medskip

\subsection{Convergence of diagonal values at local eigenvalues} 

\begin{proposition} \label{prop_203}
Assume RH.
Let $\gamma_1<\gamma_2<\cdots$ be
the distinct positive zeros in $\Gamma$,
and enumerate the positive elements of $\Gamma_{a,\pi}$ as
$0<\gamma_{a,1}<\gamma_{a,2}<\cdots$.
Then, for each fixed $n\ge1$,
\begin{equation}\label{eq_207}
 \gamma_{a,n}\searrow\gamma_n,\qquad
 J_a(\gamma_{a,n},\gamma_{a,n})
 \longrightarrow\frac1{m_{\gamma_n}}
 \quad(a\to\infty).
\end{equation}
Here $\searrow$ denotes monotone nonincreasing convergence.
For the negative eigenvalues, $-\gamma_{a,n}\nearrow-\gamma_n$, and the diagonal values have the same limit.
\end{proposition}

\begin{proof}
The monotone convergence of the eigenvalues $\gamma_{a,n}$ was proved in \cite[Theorem~1.1]{Su26b}. 
We now determine the limit of the diagonal values.
Under RH, write $\mathcal H_W$ for
the completion of $C_c^\infty(\R)$
with respect to $Q_W$.
In \cite[\S2.3]{Su26b},
a self-adjoint operator $\mathsf D_{\pi/2}$
on $\mathcal H_W$ is constructed;
its eigenvalues are precisely the points of $\Gamma$, and all are simple.
By \cite[Theorem~3.3, (4.1)]{Su26b},
if we embed the local space $\mathcal H(\mathsf A_a)$ isometrically into $\mathcal H_W$,
transport $(\overline{\mathscr D}_{a,\pi})^{-1}$
by this embedding, and extend it by zero on the orthogonal complement,
the resulting operator converges to $\mathsf D_{\pi/2}^{-1}$
in Hilbert--Schmidt norm.
The eigenvalue $1/\gamma_n$ of the latter operator is isolated and simple.
Let $P_n$ be the orthogonal projection onto the eigenspace of the global operator 
corresponding to $1/\gamma_n$, and let $P_{a,n}$ be the orthogonal projection 
onto the eigenspace of the local operator corresponding to $1/\gamma_{a,n}$. 
Convergence of the resolvents and the formula for the Riesz projections then give
$\|P_{a,n}-P_n\|\to 0$ $(a\to\infty)$.
Fix $f\in C_c^\infty(\R)$ such that
$\widehat f(\gamma_n)\ne0$.
For sufficiently large $a$, we have $f\in C_c^\infty(-a,a)$,
so \eqref{eq_202} and
Proposition~\ref{prop_202} give
\[
 \|P_{a,n}f\|_{\mathcal H_W}^2
 =\frac{|\widehat f(\gamma_{a,n})|^2}
        {J_a(\gamma_{a,n},\gamma_{a,n})}.
\]
On the other hand, as described in \cite[\S2.3]{Su26b},
$P_n$ is the orthogonal projection onto the eigenspace for $1/\gamma_n$
associated with $\gamma_n$.
Hence
\[
 \|P_nf\|_{\mathcal H_W}^2
 =m_{\gamma_n}|\widehat f(\gamma_n)|^2.
\]
It follows that
\[
 \frac{|\widehat f(\gamma_{a,n})|^2}
        {J_a(\gamma_{a,n},\gamma_{a,n})}
 \longrightarrow
 m_{\gamma_n}|\widehat f(\gamma_n)|^2.
\]
Since $\gamma_{a,n}\to\gamma_n$ and
$\widehat f(\gamma_n)\ne0$,
we obtain the second limit in \eqref{eq_207}.
By the reflection symmetry $J_a(-w,-z)=J_a(w,z)$ and \eqref{eq_203}, we have $W(a,\pi;-z)=W(a,\pi;z)$. The assertion for negative eigenvalues now follows from $J_a(-t,-t)=J_a(t,t)$ and $m_{-\gamma}=m_\gamma$.
\end{proof}

\section{RH and uniform bounds for reproducing kernels} \label{section_3}

In this section, we use the point evaluations from Section~\ref{section_2} 
and the behavior of the lowest eigenvalue $\lambda_a$ of $\mathsf A_a$ as $a$ varies.
\medskip

\noindent
{\bf Proof of Theorem~\ref{thm_3}.}~
Under RH, $Q_W^R$ is nondegenerate; 
see Yoshida~\cite[Theorem~2]{Yo92} and the discussion after Theorem~1.3 in~\cite{Su26a}. 
Moreover, $\mathsf A_R$ is nonnegative
and has compact resolvent,
so $\lambda_R>0$ for every $R>0$.
For $0<a<R$, we have $\lambda_a\geq\lambda_R>0$ 
by the monotonicity established in Section~\ref{section_2}. 
Hence, by the Cauchy--Schwarz inequality,
\[
|\widehat u(t)|^2\le 2a\|u\|_{L^2}^2
\le\frac{2R}{\lambda_R}Q_W^a(u,u)
\]
for every $u\in C_c^\infty(-a,a)$ and $t\in\R$.
Since $C_c^\infty(-a,a)$ is a core for $Q_W^a$, 
this estimate extends to every $u\in\mathcal H(\mathsf A_a)$. 
By the reproducing property of $J_a$ and the Cauchy--Schwarz inequality, we have 
\[
 J_a(t,t)=\sup_{0\ne u\in\mathcal H(\mathsf A_a)}
 \frac{|\widehat u(t)|^2}{Q_W^a(u,u)}
 \le\frac{2R}{\lambda_R}, 
\]
where the equality follows by taking $\widehat u=J_a(t,\cdot)$. 

Conversely, suppose that the uniform bound holds and RH is false. 
By \cite[Theorems~1.3 and 1.4]{Su26a}, 
$\lambda_a$ is continuous for $a>0$ and positive for sufficiently small $a>0$. 
Since RH is false, $\lambda_b<0$ for some finite $b>0$; 
see the discussion after Theorem~1.3 in~\cite{Su26a}.
Since $\lambda_a$ is monotone by extension by zero,
\[
 a_\ast:=\inf\{a>0:\lambda_a\le0\}
\]
satisfies $0<a_\ast<\infty$,
$\lambda_a>0$ $(0<a<a_\ast)$, and $\lambda_{a_\ast}=0$.
As observed in Section~\ref{section_2}, $\mathsf A_{a_\ast}$ has compact resolvent,
so there is a nonzero $u_\ast\in\ker\mathsf A_{a_\ast}$.
In particular, $Q_W^{a_\ast}(u_\ast,u_\ast)=0$.

Since $C_c^\infty(-a_\ast,a_\ast)$ is a core for $Q_W^{a_\ast}$,
we can choose $u_n\in C_c^\infty(-a_\ast,a_\ast)$ converging to $u_\ast$
in the form norm given by the square root of $Q_W^{a_\ast}(u,u)+\|u\|_{L^2}^2$.
As $Q_W^{a_\ast}$ is nonnegative and
$Q_W^{a_\ast}(u_\ast,u_\ast)=0$,
we have $
 \|u_n-u_\ast\|_{L^2}^2~\to~ 0$ and
$
 Q_W^{a_\ast}(u_n,u_n)~\to~0.$
For each $u_n$, choose $0<a_n<a_\ast$ so that its support is contained in the corresponding interval.
The forms agree under extension by zero, so
$Q_W^{a_n}(u_n,u_n)=Q_W^{a_\ast}(u_n,u_n)$.
Applying the uniform bound with $R=a_\ast$, we obtain, for every $t\in\R$,
\[
 |\widehat u_n(t)|^2
 \le J_{a_n}(t,t)Q_W^{a_n}(u_n,u_n)
 \le C_{a_\ast}Q_W^{a_\ast}(u_n,u_n)~\to~0.
\]
On the other hand, extending both functions by zero gives
\[
 \sup_{t\in\R}|\widehat u_n(t)-\widehat u_\ast(t)|
 \le\|u_n-u_\ast\|_{L^1}
 \le\sqrt{2a_\ast}\,\|u_n-u_\ast\|_{L^2}~\to~0.
\]
Thus $\widehat u_\ast(t)=0$ for every $t\in\R$,
and the injectivity of the Fourier transform gives $u_\ast=0$, a contradiction.
\hfill $\Box$
\medskip

\begin{corollary} \label{cor_301}
Let $R>0$ and suppose that $\lambda_a>0$ for every $0<a<R$.
Then $J_a(t,t)$ is uniformly bounded for $(a,t) \in (0,R)\times \R$ if and only if
$\lambda_R>0$.
\end{corollary}
\begin{proof}
If $\lambda_R>0$, the bound for point evaluations in the proof of Theorem~\ref{thm_3} gives
$J_a(t,t)\le2R/\lambda_R$.
Conversely, the continuity of $\lambda_a$ implies that $\lambda_R\ge0$.
If $\lambda_R=0$, choose a nonzero eigenfunction of $\mathsf A_R$
with eigenvalue zero and repeat the argument above with $a_\ast=R$,
using approximation by the form core and point evaluations.
This contradicts the uniform bound.
\end{proof}

\section{de Branges spaces and canonical systems}\label{section_4}

In the next section, we associate a chain of de Branges spaces with a single canonical system.
To prepare for this, we briefly recall the ordering theorem of de Branges spaces,
results on Weyl functions of canonical systems,
and the theorem of Krein--Langer on extensions of screw functions.

\subsection{de Branges spaces}\label{section_4_1}

The notation and basic facts in this subsection follow de Branges \cite{dB}.
For an entire function $F$, write $F^\sharp(z)=\overline{F(\bar z)}$.
An entire function $E$ is called a Hermite--Biehler function if
$|E^\sharp(z)|<|E(z)|$ in the upper half-plane $\C_+=\{z:\Im z>0\}$.
The de Branges space $\mathcal H(E)$ generated by $E$ consists of the entire functions $F$
such that $\int_\R |F(t)/E(t)|^2\,dt<\infty$ and both $F/E$ and $F^\sharp/E$
are of bounded type and have nonpositive mean type in $\C_+$.
Here a function is of bounded type if it is the quotient of two bounded holomorphic functions,
and the mean type of a nonzero function $f$ is $\limsup_{y\to\infty}y^{-1}\log|f(iy)|$.
The integral above gives $\|F\|_{\mathcal H(E)}^2$ \cite[\S19]{dB}.

For a Hermite--Biehler function $E$,
the functions $A=(E+E^\sharp)/2$ and $B=i(E-E^\sharp)/2$ are real entire,
and $E=A-iB$.
The de Branges space $\mathcal H(E)$ is a reproducing kernel Hilbert space.
With the inner product linear in the first variable, its reproducing kernel, satisfying
$F(w)=\langle F,K_E(w,\cdot)\rangle$, is
\begin{equation} \label{eq_401}
 K_E(w,z)=
 \frac{\overline{A(w)}B(z)-\overline{B(w)}A(z)}{\pi(z-\bar w)}
 =\frac{\overline{E(w)}E(z)-\overline{E^\sharp(w)}E^\sharp(z)}
 {2\pi i(\bar w-z)}
\end{equation}
\cite[Theorem~19]{dB}. Fixing the point $i$, we can choose a generating function from the reproducing kernel.
Indeed, set $F(z)=(z+i)K_E(i,z)\sqrt{\pi/K_E(i,i)}$ and write it as
$cE+dE^\sharp$.
By \eqref{eq_401}, we have $|E(i)|^2-|E^\sharp(i)|^2=4\pi K_E(i,i)$,
so $|c|^2-|d|^2=1$. Thus
$|F(z)|^2-|F^\sharp(z)|^2
=|E(z)|^2-|E^\sharp(z)|^2$,
and since $E$ is a Hermite--Biehler function,
so is $F$.
Substituting this relation into the numerator in \eqref{eq_401} gives
$K_F=K_E$.
Hence $F$ generates the same de Branges space isometrically.

A de Branges subspace is a closed Hilbert subspace that is itself a de Branges space
with the induced norm.
The following ordering theorem describes the inclusion relations 
between de Branges spaces \cite[Theorem~35]{dB}.

\begin{theorem}[de Branges ordering theorem]\label{thm_401}
Suppose that two de Branges spaces $\mathcal H(E_1),\mathcal H(E_2)$ embed isometrically,
by restriction to the real axis, into the same $L^2$ space for a positive measure.
If $E_2/E_1$ is of bounded type in the upper half-plane and has neither zeros nor poles
on the real axis, then one space is contained in the other.
\end{theorem}

\subsection{Canonical systems}
\label{section_4_2}

The notation and basic facts in this subsection follow
Krein--Langer~\cite{KL} and Winkler~\cite{Win15}.
A Hamiltonian $H$ is a locally integrable, real symmetric, nonnegative
$2\times2$ matrix-valued function
that is not identically $0$ on any open interval.
Let $J$ be the matrix defined in \eqref{eq_105}.
The system of first-order linear ordinary differential equations
\[
JY'=-zHY 
\]
parametrized by $z \in \C$, is called a canonical system with Hamiltonian $H$.
Define the fundamental matrix $U(x,z)$ by
\begin{equation} \label{eq_402}
 \partial_x U(x,z)J=zU(x,z)H(x),
 \qquad U(0,z)=I_2.
\end{equation}
The transpose of each row solves the canonical system.
The entries of $U(x,z)$ are real entire functions of $z$,
and $\det U(x,z)=1$.
The fundamental matrix on $[x,y]$, normalized to the identity matrix $I_2$
at the left endpoint, is $U(x,z)^{-1}U(y,z)$.
Writing the second row of the fundamental matrix as $(-B_x(z),A_x(z))$,
we obtain from \eqref{eq_402} that
\[
 \frac{\overline{A_x(w)}B_x(z)
       -\overline{B_x(w)}A_x(z)}{z-\bar w}
 =
 \int_0^x (-B_t(z),A_t(z))H(t)
 \binom{-\overline{B_t(w)}}{\overline{A_t(w)}}\,dt.
\]
If $\int_0^x h_{22}(t)\,dt>0$, then
$A_x-iB_x$ is a Hermite--Biehler function.
Under the convention in \eqref{eq_401}, the expression above is the reproducing kernel
of the de Branges space generated by
$\sqrt\pi(A_x-iB_x)$
\cite[(22.9), (22.25)--(22.30)]{Win15}. 
If $\operatorname{tr}H(x)=1$ a.e.,
we call $H$ a \emph{trace-normed Hamiltonian}.

\subsection{Weyl functions and boundary conditions at the right endpoint}
The notation and basic facts in this subsection follow
Winkler~\cite{Win95, Win15}.
A Nevanlinna function is holomorphic in the upper half-plane, has nonnegative imaginary part,
and is extended to the lower half-plane by conjugate reflection.
Every Nevanlinna function $q\ne\infty$ has a unique representation
\[
 q(z)=c+dz+\int_{-\infty}^{\infty}
 \left(\frac{1}{\lambda-z}
       -\frac{\lambda}{1+\lambda^2}\right)d\nu(\lambda),
 \quad c\in\R,\quad d\ge0,
 \quad
 \int_{-\infty}^{\infty}
 \frac{d\nu(\lambda)}{1+\lambda^2}<\infty.
\]
Here $\nu$ is a nonnegative Borel measure,
called the spectral measure of $q$ 
\cite[(22.14), (22.15)]{Win15}.
We write the action of a matrix $U=(U_{jk})$ by a linear fractional transformation as
\[
 U\cdot\Omega
 =
 \frac{U_{11}\Omega+U_{12}}
      {U_{21}\Omega+U_{22}}.
\]
For $\Omega=\infty$, this means $U_{11}/U_{21}$;
if the denominator vanishes, the value is $\infty$.

A canonical system on the half-line is said to be in the \emph{limit-point case}
at infinity if, for each $z\in\C_+$, the limit
\[
 \lim_{x\to\infty}U(x,z)\cdot\Omega(z)
\]
exists and is independent of the choice of Nevanlinna function $\Omega$
or $\Omega=\infty$.
The function $q$ defined by this limit is called the \emph{Weyl function}.
It is either a Nevanlinna function or $\infty$.

\begin{proposition} \label{prop_402}
Every canonical system associated with a trace-normed Hamiltonian
$H$ on the half-line $[0,\infty)$ is in the limit-point case at infinity.
Conversely, every Nevanlinna function $q$, as well as
$q=\infty$, is the Weyl function of a trace-normed Hamiltonian $H$
on the half-line, unique up to equality a.e.
The function $q=\infty$ corresponds to
$H=\operatorname{diag}(1,0)$.
\end{proposition}

The first assertion follows from \cite[Theorem~2.1]{Win95},
and the second from \cite[Theorem~1]{Win95}.
By uniqueness in the inverse problem,
$H$ is diagonal whenever its Weyl function $q$ satisfies $q(-z)=-q(z)$.
Indeed, conjugating $H$ by $\operatorname{diag}(1,-1)$
gives the fundamental matrix obtained by conjugating $U(x,-z)$ by $\operatorname{diag}(1,-1)$,
whose Weyl function is $-q(-z)$.
Thus, again by uniqueness, $H$ is invariant under this conjugation,
so its off-diagonal entries vanish a.e.

Consider an extension of a Hamiltonian $H$ from a finite interval $[0,x]$
to $[0,\infty)$.
Let $\Omega$ be the Weyl function of the canonical system
on the added interval $[x,\infty)$,
defined by normalizing its fundamental matrix to
$I_2$ at the left endpoint $x$.

\begin{proposition} \label{prop_403}
The Weyl function of the canonical system on the half-line obtained by this extension
is $U(x,z)\cdot\Omega(z)$.
\end{proposition}

This follows from \cite[Lemma~2.2]{Win95}.
The \emph{spectral measure} of a canonical system on the half-line is the spectral measure
of its Weyl function.
The spectral measures of a canonical system on a finite interval $[0,L]$
are those of $U(L,z)\cdot\Omega(z)$ as $\Omega$ ranges over
all Nevanlinna functions and $\Omega=\infty$,
excluding the case $U(L,z)\cdot\Omega(z)\equiv\infty$.

\subsection{Screw functions and parameters for extensions}
\label{section_4_4}

This subsection follows
Krein--Langer \cite[\S5]{KL}.
Let $0<a\le\infty$. 
A continuous function $g$ on $(-2a,2a)$
is called a screw function on $(-2a,2a)$ if it satisfies
$g(-t)=\overline{g(t)}$ and the kernel
\[
 G_g(s,t)=g(t-s)-g(t)-g(-s)+g(0),
 \qquad |s|,|t|<a
\]
is positive semidefinite
\cite[(5.1)]{KL}.
Here positive semidefiniteness means that the matrix
$(G_g(t_j,t_k))_{j,k}$ is positive semidefinite for every finite collection of points.
If $0<a<\infty$, then $g$ extends to a screw function on the whole real line
\cite[Corollary~5.2]{KL}.
Let $C_0(-a,a)$ denote the space of continuous functions
with compact support in $(-a,a)$, and set
\[
 L(G_g):=\left\{\varphi\in C_0(-a,a):
                  \int_{-a}^a\varphi(t)\,dt=0\right\}.
\]
Equip this space with the semidefinite inner product
\[
 \langle\varphi,\psi\rangle_{G_g}
 :=\int_{-a}^a\!\int_{-a}^a
 G_g(s,t)\varphi(t)\overline{\psi(s)}\,dt\,ds.
\]
Let $\mathcal H(G_g)$ be the Hilbert space obtained by taking the quotient
by the subspace of functions with zero seminorm and completing.

By \cite[Theorem~5.1]{KL}, every extension $\widetilde g$ of a screw function $g$
on a finite interval $(-2a,2a)$ to a screw function on the whole real line
has the representation
\[
 \widetilde g(t)
 =\widetilde g(0)+i\gamma t
 +\int_{-\infty}^{\infty}
 \left(e^{it\lambda}-1
       -\frac{it\lambda}{1+\lambda^2}\right)
       \frac{d\nu(\lambda)}{\lambda^2},
 \qquad \gamma\in\R,
\]
where $\nu$ is a nonnegative Borel measure satisfying
$\int_{-\infty}^{\infty}(1+\lambda^2)^{-1}d\nu(\lambda)<\infty$,
and the integrand at $\lambda=0$ is defined to be
$-t^2/2$.
The measure $\nu$ is called the
\emph{spectral measure} of $\widetilde g$.
By \cite[(5.16), (5.17)]{KL}, the transform
\begin{equation} \label{eq_403}
 \Phi_1(\varphi;\lambda)
 =\int_{-a}^a
   \frac{e^{it\lambda}-1}{\lambda}\varphi(t)\,dt,
 \qquad
 \left.\frac{e^{it\lambda}-1}{\lambda}\right|_{\lambda=0}
 =it
\end{equation}
satisfies
\[
 \langle\varphi,\psi\rangle_{G_g}
 =\int_{-\infty}^{\infty}\Phi_1(\varphi;\lambda)
       \overline{\Phi_1(\psi;\lambda)}\,d\nu(\lambda).
\]
Hence $\Phi_1$ extends to an isometry
from $\mathcal H(G_g)$ into $L^2(\R,d\nu)$.

From now on, assume that $g(0)=0$.
For an extension $\widetilde g$ of $g$ to a screw function
on the whole real line, set
\[
 q_{\widetilde g}(z)
 =-iz^2\int_0^\infty e^{-izt}\widetilde g(t)\,dt,
 \qquad \Im z<0.
\]
The Nevanlinna function obtained by conjugate reflection has the representation
\[
 q_{\widetilde g}(z)
 =c+\int_{-\infty}^{\infty}
 \left(\frac{1}{\lambda-z}
       -\frac{\lambda}{1+\lambda^2}\right)d\nu(\lambda),
 \quad c\in\R,\quad
 \int_{-\infty}^{\infty}\frac{d\nu(\lambda)}{1+\lambda^2}<\infty
\]
\cite[(5.5)--(5.7)]{KL}.
In particular, if $\widetilde g$ is real and even,
then $c=0$, the measure $\nu$ is symmetric, and
\[
 q_{\widetilde g}(z)
 =z\int_{-\infty}^{\infty}\frac{d\nu(\lambda)}{\lambda^2-z^2}.
\]
We have the following result of Krein and Langer. 

\begin{theorem} \label{thm_404}
Let $0<a<\infty$ and $g(0)=0$, and suppose that
a screw function $g$ on $(-2a,2a)$ has at least two distinct extensions
to screw functions on the whole real line.
Then there is a unique trace-normed Hamiltonian $H=(h_{jk})$
on a finite interval $[0,L]$ such that, for its fundamental matrix $U(x,z)$,
the formula
\[
 q_{\widetilde g}(z)=U(L,z)\cdot\Omega(z)
\]
gives a one-to-one correspondence between all extensions $\widetilde g$
and all Nevanlinna functions $\Omega$ together with $\Omega=\infty$.
The spectral measures of this canonical system are precisely
the spectral measures of extensions of $g$.
Moreover,
\[
 \int_0^x h_{22}(t)\,dt>0\quad(0<x\le L),\qquad
 \int_x^L\sqrt{\det H(t)}\,dt>0\quad(0\le x<L),
\]
\[
 a=\int_0^L\sqrt{\det H(t)}\,dt.
\]
\end{theorem}

This follows from \cite[Theorems~5.4--5.6]{KL}.
The matrix $U(x,z)$ in the theorem follows the convention in \eqref{eq_402},
with $U(0,z)=I_2$.
Its value $U(L,z)$ at the right endpoint is a $2\times2$ matrix of entire functions
that represents the Nevanlinna functions corresponding to extensions
by linear fractional transformations.
It is called the \emph{resolvent matrix} of $g$
\cite[Theorem~5.4]{KL}.

\medskip
We next define orthogonal spectral measures.
The operator
$\varphi\mapsto i\varphi'$
on $L(G_g)$
(with $\varphi,\varphi'\in L(G_g)$)
gives rise to a densely defined symmetric operator on $\mathcal H(G_g)$.
When a self-adjoint extension is realized in $\mathcal H(G_g)$ itself,
the corresponding extension $\widetilde g$ and spectral measure $\nu$
are called \emph{orthogonal} \cite[\S5.3]{KL}.
Orthogonal spectral measures admit 
the following characterization \cite[\S5.3, (5.16)--(5.18)]{KL}.

\begin{proposition} \label{prop_405}
Under the assumptions of Theorem~\ref{thm_404}, the following are equivalent:
the extension $\widetilde g$ (or its spectral measure $\nu$) is orthogonal;
the isometry $\mathcal H(G_g)\to L^2(\R,d\nu)$ in \eqref{eq_403}
is surjective; and the parameter $\Omega$ in Theorem~\ref{thm_404}
is a real constant or $\infty$.
\end{proposition}

\section{Chains of de Branges subspaces and a canonical system}\label{section_5}

Throughout this section, fix $R>0$ and assume that $\lambda_R>0$.
Extension by zero preserves the value of the local Weil quadratic form, so
$\lambda_a\ge\lambda_R>0$ for $0<a\le R$.
We apply the results of Section~\ref{section_4} to $\mathcal B_a=\widehat{\Ha}$
to construct a canonical system realizing this chain and determine the point corresponding to each $\mathcal B_a$.

\subsection{Chains of de Branges spaces and generating functions}

For every $0<a\le R$, we can argue as in \cite[\S7.9]{Su26a}
to show that $\mathcal B_a$ is a de Branges space. 
Indeed, the equality of norms used there under RH follows here from
\eqref{eq_205}, which holds under $\lambda_a>0$.
Moreover, $J_a(z,z)>0$ implies that the space has no common zeros.

\begin{proposition} \label{prop_501}
Set $\mathcal B_0=\{0\}$.
For $0\le a<b\le R$, the space $\mathcal B_a$ is a proper closed subspace
of $\mathcal B_b$, with the same norm.
Moreover, in the norm of $\mathcal B_R$,
\begin{equation} \label{eq_501}
 \overline{\bigcup_{b<a}\mathcal B_b}
 =\mathcal B_a\quad(0<a\le R),
 \qquad
 \bigcap_{a<b\le R}\mathcal B_b
 =\mathcal B_a\quad(0\le a<R).
\end{equation}
Furthermore, every nonzero de Branges subspace of $\mathcal B_R$
with no common real zeros equals $\mathcal B_a$
for a unique $a\in(0,R]$.
\end{proposition}
\begin{proof}
We first prove the strict inclusion for $0\le a<b\le R$.
Since $C_c^\infty(-a,a)\subset C_c^\infty(-b,b)$ and the Weil quadratic forms
agree on the smaller space, $\mathcal B_a$ is a closed subspace of $\mathcal B_b$
with the same norm.
The inclusion is strict:
the Fourier transform of a nonzero test function supported in $(a,b)$
belongs to $\mathcal B_b$ but not to $\mathcal B_a$,
since the inverse Fourier transform of every element of $\mathcal B_a$
is supported in $[-a,a]$.

We next prove \eqref{eq_501}.
The first identity follows immediately from the definition of $\mathcal B_a$,
since every element of $C_c^\infty(-a,a)$ belongs to $C_c^\infty(-b,b)$
for some $b<a$.
For the second identity, the right-hand side is contained in the left-hand side
by the inclusion just proved.
Conversely, the inverse Fourier transform of an element of the left-hand side
is supported in $[-b,b]$ for every $b>a$, and hence in $[-a,a]$.
By \eqref{eq_201}, this inverse Fourier transform $u$ also satisfies
$\int_{-\infty}^\infty\log(2+|t|)|\widehat u(t)|^2\,dt<\infty$.
If $a>0$, a function satisfying this condition and supported in $[-a,a]$
can be approximated by elements of $C_c^\infty(-a,a)$
by slightly contracting the interval and then smoothing.
These operations are continuous in the $H^{\log}$ norm. 
Since this norm is equivalent to the norm induced by $Q_W^a$,
we obtain $\widehat u\in\mathcal B_a$.
When $a=0$, the inverse Fourier transform is supported
at the single point $0$ and hence vanishes as an $L^2$ function.
Thus the intersection is $\{0\}=\mathcal B_0$.
This proves \eqref{eq_501}.
The orthogonal projections onto an increasing or decreasing family of closed subspaces
converge strongly. It follows that the orthogonal projection from $\mathcal B_R$
onto $\mathcal B_a$ is strongly continuous for $0\le a\le R$.

Finally, let $\mathcal X$ be a nonzero de Branges subspace of $\mathcal B_R$
with no common real zeros.
We show that $\mathcal X$ and every $\mathcal B_b$ are comparable by inclusion.
If $E$ generates $\mathcal B_R$, both spaces embed isometrically into
$L^2(\R,|E(t)|^{-2}\,dt)$.
Setting $w=i$ in \eqref{eq_401}, we can choose a generating function
for each space to be a constant multiple of $(z+i)K(i,z)$,
where $K$ is its reproducing kernel (cf. Section~\ref{section_4_1}).
The function $K(i,\cdot)\in\mathcal B_R$ is the Fourier transform of an $L^2$ function
supported in $[-R,R]$, so $e^{iRz}K(i,z)$ is bounded in the upper half-plane.
Thus these generating functions are of bounded type in the upper half-plane.
They have no zeros on the real axis, so their quotient is also of bounded type
and has neither zeros nor poles on the real axis.
Theorem~\ref{thm_401} therefore shows that the two de Branges spaces are comparable.

Set
\[
 a=\sup\bigl(\{0\}\cup
 \{b\in(0,R]:\mathcal B_b\subset\mathcal X\}\bigr).
\]
By the definition of $a$ and the first identity in \eqref{eq_501}, 
we have $\mathcal B_a\subset\mathcal X$.
If $a<b\le R$, then $\mathcal B_b\subset\mathcal X$ 
would contradict the definition of $a$.
Since $\mathcal X$ and $\mathcal B_b$ are comparable by inclusion, 
we must have $\mathcal X\subset\mathcal B_b$. 
For $a<R$, the second identity in \eqref{eq_501} gives
$\mathcal X\subset\mathcal B_a$.
For $a=R$, we also have $\mathcal X\subset\mathcal B_R
=\mathcal B_a$.
Hence $\mathcal X=\mathcal B_a$.
Since $\mathcal X$ is nonzero, $a>0$.
The strict inclusion proved at the beginning shows that this $a$ is unique.
\end{proof}

\begin{proposition} \label{prop_502}
Let $0<a\le R$ and define
\[
 E_a(z):=\sqrt{\frac{\pi}{J_a(i,i)}}\,(z+i)J_a(i,z).
\]
Then $E_a$ is a Hermite--Biehler function generating $\mathcal B_a$.
Using \eqref{eq_203}, define
\begin{equation} \label{eq_502}
 A_a(z):=\frac{W(a,\pi;z)}{W(a,\pi;0)},\qquad
 B_a(z):=zJ_a(0,z).
\end{equation}
These are respectively even and odd real entire functions, with $A_a(0)=1$ and $B_a(0)=0$.
Moreover,
\[
\mathcal E_a:= \sqrt\pi(A_a-iB_a)
\]
also generates $\mathcal B_a$, whose reproducing kernel is
\begin{equation} \label{eq_503}
 J_a(w,z)=
 \frac{\overline{A_a(w)}B_a(z)-\overline{B_a(w)}A_a(z)}
      {z-\bar w}.
\end{equation}
The pair of real entire functions $(A_a,B_a)$ satisfying \eqref{eq_503}
is uniquely determined by the requirements that $A_a$ be even, $B_a$ be odd,
and $A_a(0)=1$, $B_a(0)=0$.
\end{proposition}
\begin{proof}
Taking $w=i$ and $K_E=J_a$ in \eqref{eq_401} shows that
$E_a$ generates $\mathcal B_a$.
Since $\mathcal B_a$ has no common real zeros, $E_a$ has no real zeros.
The reflection $F(z)\mapsto F(-z)$ is isometric by the symmetry of \eqref{eq_102}.
Together with the isometry $F\mapsto F^\sharp$, this gives
$J_a(i,-z)=J_a(-i,z)=\overline{J_a(i,\bar z)}$, and hence $E_a(-z)=-E_a^\sharp(z)$.
In particular, $E_a(iy)$ is purely imaginary for real $y$.
Since $E_a$ is a Hermite--Biehler function without real zeros,
$E_a(iy)\ne0$ for $y\ge0$.
By continuity and the equality 
$E_a(i)=2i\sqrt{\pi J_a(i,i)}$,
the imaginary part of $E_a(0)$ is also positive.
Thus $J_a(i,0)>0$.
Set
\[
 \alpha_a=J_a(i,0),\qquad \beta_a=\frac{J_a(i,i)}{J_a(i,0)}.
\]
Substituting $\theta=\pi$ into \eqref{eq_203} and applying
\eqref{eq_401} with $E=E_a$ and $w=0$, we obtain
\[
 W(a,\pi;0)=-2i\alpha_a,\qquad
 \sqrt{\frac{J_a(i,i)}{\pi}}\,
 \frac{E_a-E_a^\sharp}{2i}=\alpha_a A_a,\qquad
 \sqrt{\frac{J_a(i,i)}{\pi}}\,
 \frac{E_a+E_a^\sharp}{2}=\beta_a B_a .
\]
Hence $A_a,B_a$ have the asserted parity and values at the origin, and
\[
 E_a=\sqrt{\frac{\pi}{J_a(i,i)}}\,(i\alpha_a A_a+\beta_a B_a).
\]
In particular, $B_a(i)=i\alpha_a$ and $A_a(i)=\beta_a$.
Since $\alpha_a\beta_a=J_a(i,i)$, substituting this formula into the expression for the reproducing kernel
gives \eqref{eq_503}.
As $J_a(w,w)>0$ in the upper half-plane,
$\mathcal E_a=\sqrt\pi(A_a-iB_a)$ is also a Hermite--Biehler function and generates the same space.

For uniqueness, evaluating \eqref{eq_503} at $w=0$ determines $B_a$.
The same equality shows that any two candidates for $A_a$ differ by a constant multiple of $B_a$.
Since $A_a$ is even and $B_a$ is odd, this difference vanishes.
\end{proof}

\subsection{Screw functions and local extensions} \label{section_5_2}

Unless otherwise stated, $g$ denotes the real even function on the real axis
associated with the Riemann zeta function in \cite[(1.8)]{Su23}.
By \cite[\S3.4]{Su23}, if
$u_t=2^{-1/2}\mathbf1_{[-|t|/2,\,|t|/2]}$, then
$g(t)=-W(u_t*\widetilde u_t)$.
In particular, $g(0)=0$.
The following result relates local restrictions of $g$
to the extension theory in Section~\ref{section_4_4}.

\begin{proposition} \label{prop_503}
Assume that $\lambda_R>0$ and let $0<a\le R$.
Then $g|_{(-2a,2a)}$ is a screw function.
Moreover, for every $\theta\in[0,2\pi)$,
\begin{equation} \label{eq_504}
 g(t)=\int_{-\infty}^\infty
 \frac{\cos(t\lambda)-1}{\lambda^2}
 \,d\mu_{a,\theta}(\lambda)
 \qquad(|t|\le2a).
\end{equation}
Here $\mu_{a,\theta}$ is the spectral measure defined in \eqref{eq_204},
and the integrand at $\lambda=0$ is defined to be $-t^2/2$.
Furthermore, the measures for the phases $0$ and $\pi$ give two distinct extensions
of $g|_{(-2a,2a)}$ to screw functions on the whole real line.
\end{proposition}
\begin{proof}
Taking the Fourier transform in the definition gives
$\widehat u_t(\lambda)=\sqrt2\sin(|t|\lambda/2)/\lambda$.
By \eqref{eq_201}, we have $u_t\in H^{\log}(-a,a)$ whenever $|t|\le2a$.
Lemma~\ref{lem_201} and the local Parseval identity
\eqref{eq_205} now give \eqref{eq_504}.

By \eqref{eq_203} and Proposition~\ref{prop_502},
$W(a,0;z)=2\beta_a B_a(z)$ and
$W(a,\pi;z)=-2i\alpha_a A_a(z)$.
The parity of $A_a,B_a$ and the identity $J_a(-t,-t)=J_a(t,t)$ show that
both $\mu_{a,0}$ and $\mu_{a,\pi}$ are symmetric measures.
Applying the estimate for the growth of the measure in the proof of Theorem~\ref{thm_1}
to each interval $2^n<|\lambda|\le 2^{n+1}$ gives
$\int_{-\infty}^\infty(1+\lambda^2)^{-1}\,
d\mu_{a,\theta}(\lambda)<\infty$.
We may therefore use $\mu_{a,\theta}$ to define a function $g_\theta(t)$
on the whole real line by the right-hand side of \eqref{eq_504}.
For both $\theta=0,\pi$, the integral converges absolutely and uniformly
on compact sets of $t$, so $g_\theta$ is continuous.
By symmetry of the measure,
\[
g_\theta(t-u)-g_\theta(t)-g_\theta(-u)+g_\theta(0)
=\int_{-\infty}^\infty
 \frac{(e^{it\lambda}-1)(e^{-iu\lambda}-1)}{\lambda^2}
 \,d\mu_{a,\theta}(\lambda).
\]
The right-hand side is the inner product in $L^2(\mu_{a,\theta})$
of the functions
$\lambda\mapsto(e^{it\lambda}-1)/\lambda$ and
$\lambda\mapsto(e^{iu\lambda}-1)/\lambda$,
with values $it$ and $iu$ at the origin, respectively.
Hence this kernel is positive semidefinite. 
Thus $g_0$ and $g_\pi$ are both screw functions on the whole real line,
and \eqref{eq_504} shows that they agree with $g$ on $(-2a,2a)$.
In particular, $g|_{(-2a,2a)}$ is a screw function.

Finally, of the two measures, only $\mu_{a,0}$
has positive mass at the origin.
This mass contributes $-t^2/(2J_a(0,0))$ to $g_0(t)$. 
Since the spectrum is discrete and has no finite accumulation points,
the contribution away from the origin is uniformly bounded by
$2\int_{\lambda\ne0}\lambda^{-2}\,d\mu_{a,\theta}$.
Hence $g_\pi$ is bounded, whereas $g_0$ is unbounded,
so the two extensions are distinct.
\end{proof}

\subsection{A common canonical system and the points $L_a$}

By Proposition~\ref{prop_503}, for each $0<a\le R$,
the function $g|_{(-2a,2a)}$ is a screw function with $g(0)=0$
and has two distinct extensions to the whole real line.
Thus Theorem~\ref{thm_404} applies.
Denote the corresponding trace-normed Hamiltonian on the finite interval $[0,L_a]$ by
\[
H_{\rm tr}=(h_{jk})_{j,k=1}^2,
\]
and the value of its fundamental matrix at $x=L_a$ by
\[
U_a(z)=(U_{a,jk}(z))_{j,k=1}^2.
\]
At this stage, the Hamiltonian is defined separately for each $a$.
By Theorem~\ref{thm_404},
\[
 \int_0^x h_{22}(y)\,dy>0\quad(0<x\le L_a),
 \qquad
 \int_0^{L_a}\sqrt{\det H_{\rm tr}(x)}\,dx=a,
\]
and $L_a$ is the first point at which the latter integral reaches $a$.
Moreover, all extensions of $g|_{(-2a,2a)}$ to screw functions on the whole real line
are parametrized by Nevanlinna functions $\Omega$ and $\Omega=\infty$.
The corresponding Weyl functions are
\begin{equation} \label{eq_505}
 q_\Omega(z)=
 \frac{U_{a,11}(z)\Omega(z)+U_{a,12}(z)}
      {U_{a,21}(z)\Omega(z)+U_{a,22}(z)},
 \qquad
 q_\infty(z)=\frac{U_{a,11}(z)}{U_{a,21}(z)}.
\end{equation}

\begin{proposition} \label{prop_504}
Assume that $\lambda_R>0$. With the notation above,
all matrices $U_a$ for $0<a\le R$ can be written as
$U_a(z)=U(L_a,z)$, where $U(x,z)$ is the fundamental matrix
of a single canonical system with a diagonal trace-normed Hamiltonian.
The restriction of this system to $[0,L_R]$ is unique under the normalization above,
and $L_a$ is strictly increasing in $a$.
Under RH, the same system can be used for all $a>0$.
\end{proposition}

\begin{proof}
We first show that $\mathcal H(G_g)$ and $\Ha$ are isomorphic.
For $u\in C_c^\infty(-a,a)$, we have $iu'\in L(G_g)$.
Integration by parts in \eqref{eq_403} gives
\[
 \Phi_1(iu';\lambda)
 =\int_{-a}^a
   \frac{e^{it\lambda}-1}{\lambda}iu'(t)\,dt
 =\widehat u(\lambda).
\]
Both sides are continuous at $\lambda=0$, so the identity also holds there.
For $u,f\in C_c^\infty(-a,a)$,
Proposition~\ref{prop_503} and \eqref{eq_403} give
\[
\langle iu',if'\rangle_{G_g}
=\int_{-\infty}^{\infty}
  \widehat u(\lambda)\overline{\widehat f(\lambda)}
  \,d\mu_{a,\pi}(\lambda)=Q_W^a(u,f),
\]
where the last equality follows from \eqref{eq_205}.
Thus $u\mapsto iu'$ is isometric.
Every smooth compactly supported function with integral zero can be written as $iu'$.
Moreover, every continuous compactly supported function with integral zero
can be approximated by smooth compactly supported functions while preserving integral zero.
Since the kernel is bounded on a common compact set,
this approximation also holds in the norm of $\mathcal H(G_g)$.
Thus $u\mapsto iu'$ extends to an isometric isomorphism
from $\Ha$ onto $\mathcal H(G_g)$.

For every $\theta$, \eqref{eq_205} gives the Parseval identity
for the transform $\Phi_1$ in \eqref{eq_403}. 
By \cite[\S5.3]{KL}, this identity shows that $\mu_{a,\theta}$
is the spectral measure of an extension of $g|_{(-2a,2a)}$.
Furthermore, Proposition~\ref{prop_202} shows that
$u\mapsto\widehat u$ is a surjective isometry
from $\Ha$ onto $L^2(\mu_{a,\theta})$.
Through the isometric isomorphism established above, $\Phi_1$ is also surjective.
Proposition~\ref{prop_405} therefore implies that
$\mu_{a,\theta}$ is an orthogonal spectral measure,
and that the corresponding $\Omega$ is a real constant or $\infty$.

We next show that the canonical systems obtained for each $a$
are restrictions of a single canonical system.
By Proposition~\ref{prop_503},
$\mu_{R,\pi}$ defines an extension of $g|_{(-2R,2R)}$
to a screw function on the whole real line.
This is a common extension of $g|_{(-2a,2a)}$ for all $0<a\le R$.
It is real and even; denote its Weyl function by $q$.
For each $a\in(0,R]$, choose the parameter $\Omega$ in \eqref{eq_505}
corresponding to this extension.
By Proposition~\ref{prop_402}, there is a trace-normed Hamiltonian
on the half-line with Weyl function $\Omega$.
For $\Omega=\infty$, use $H=\operatorname{diag}(1,0)$.
Attaching this system to the right of the canonical system on $[0,L_a]$
gives a system with Weyl function $q$, by Proposition~\ref{prop_403}.
By Proposition~\ref{prop_402}, the trace-normed Hamiltonian
on the half-line corresponding to $q$ is unique.
Thus the canonical systems obtained on $[0,L_a]$ for each $a$
are all restrictions of the same canonical system on the half-line.
Henceforth, let $H_{\rm tr}$ denote this common Hamiltonian,
and write $U(x,z)$ for its fundamental matrix.
Then
\[
 U_a(z)=U(L_a,z).
\]
For this Hamiltonian $H_{\rm tr}$, we have
\[
 \int_0^{L_a}\sqrt{\det H_{\rm tr}(x)}\,dx=a,
\]
so $a<b$ implies $L_a<L_b$.
Moreover, since $\mu_{R,\pi}$ is symmetric, the spectral representation
in Section~\ref{section_4_4} gives
\[
 q(z)=z\int_{-\infty}^{\infty}
       \frac{d\mu_{R,\pi}(\lambda)}{\lambda^2-z^2}.
\]
Hence $q$ is odd.
By the argument immediately following Proposition~\ref{prop_402}, 
the Hamiltonian $H_{\rm tr}$ is diagonal.
\medskip

If $R$ is changed, the systems agree on their common interval,
since the Hamiltonian in Theorem~\ref{thm_404} is unique
for the same $g|_{(-2a,2a)}$.
Under RH, we have $\lambda_R>0$ for every $R>0$.
Since the Hamiltonian is diagonal and trace-normed,
$\sqrt{\det H_{\rm tr}}\le1/2$, and hence $L_a\ge2a$.
Thus $L_a\to\infty$ as $a\to\infty$,
and we obtain a single system on the entire half-line.
\end{proof}

We have shown that the canonical systems constructed for each $a$
are restrictions of one common canonical system.
We next identify the boundary conditions corresponding to $\mu_{a,0}$ and $\mu_{a,\pi}$,
and show that the second row of its fundamental matrix
is expressed in terms of $A_a,B_a$ from Proposition~\ref{prop_502}.

\begin{proposition} \label{prop_505}
For the common canonical system obtained in Proposition~\ref{prop_504},
we have, for every $0<a\le R$,
\[
 (U_{a,21}(z),U_{a,22}(z))=(-B_a(z),A_a(z)).
\]
Moreover, under the condition at the left endpoint in
Theorem~\ref{thm_404}, the full fundamental matrix $U_a$
is uniquely determined.
\end{proposition}
\begin{proof}
Throughout this proof, write $U_{jk}=U_{a,jk}(z)$.
By Proposition~\ref{prop_504}, the common Hamiltonian $H_{\rm tr}$ is diagonal,
and the fundamental matrix satisfies $U(0,z)=I_2$.
It follows from \eqref{eq_402} that
$U_{11},U_{22}$ are even and $U_{12},U_{21}$ are odd.
For $\Omega=r\in\R\cup\{\infty\}$, \eqref{eq_505} gives
$-q_r(-z)=q_{-r}(z)$, with $-\infty=\infty$.
If $q_r$ corresponds to an orthogonal extension that is real and even, then $q_r$ is odd.
Since the linear fractional transformation is injective, $r=-r$,
so $r=0$ or $r=\infty$.
We have $U_a(0)=I_2$, and the integral equation associated with
\eqref{eq_402} gives
\[
 \partial_zU_{a,21}(0)
 =-\int_0^{L_a}h_{22}(x)\,dx<0.
\]
Thus $q_\infty=U_{11}/U_{21}$ has a pole at the origin,
whereas $q_0=U_{12}/U_{22}$ is holomorphic there.
Of the two symmetric measures in Proposition~\ref{prop_503},
only $\mu_{a,0}$ has positive mass at the origin.
Hence phase $0$ corresponds to $r=\infty$,
and phase $\pi$ to $r=0$.

Since $\det U_a=1$, the numerator and denominator of each of
$q_\infty$ and $q_0$ have no common zeros.
The poles of these Nevanlinna functions occur precisely at the points 
where the corresponding spectral measures have positive mass. 
They are all simple, with residues equal to the negatives of the corresponding masses.
Thus the zeros of $U_{21}$ and $U_{22}$ are respectively
$\Gamma_{a,0}$ and $\Gamma_{a,\pi}$, and are all simple.
On the other hand, \eqref{eq_203} and Proposition~\ref{prop_502} show that
the zeros of $B_a,A_a$ are also respectively
$\Gamma_{a,0},\Gamma_{a,\pi}$.
By \eqref{eq_503} and $J_a(t,t)>0$,
these zeros are simple as well.

Successive substitution in the integral equation for the fundamental matrix
bounds $U_{jk}$ by a constant multiple of
\[
\exp\left(|z|\int_0^{L_a}\|H_{\rm tr}(x)\|\,dx\right).
\]
Thus $U_{21},U_{22}$ are entire functions of finite exponential type.
The functions $A_a,B_a$ are also of finite exponential type,
by their representations in terms of Fourier transforms of compactly supported functions.
Hadamard factorization therefore shows that the ratio of each pair of corresponding functions 
is the exponential of a polynomial of degree at most one. 
Since $U_{21}$ and $B_a$ are odd, while $U_{22}$ and $A_a$ are even, 
both ratios are even functions and hence constant.

Since $U_{a,22}(0)=A_a(0)=1$, we have $U_{a,22}=A_a$.
Comparing residues of $q_\infty$ at the origin gives
\[
 -\frac{1}{\partial_zU_{a,21}(0)}
 =\mu_{a,0}(\{0\})
 =\frac{1}{J_a(0,0)}
 =\frac{1}{B_a'(0)}.
\]
Since both $U_{21}$ and $B_a$ are odd,
we obtain $U_{a,21}=-B_a$.

Finally, we verify uniqueness of the first row.
By \cite[Lemma~1.2]{KL}, among normalized resolvent matrices
with a prescribed second row, there is a unique one with no left factor of the form
$\left(\begin{smallmatrix}1&dz\\0&1\end{smallmatrix}\right)$
with $d>0$.
Such a left factor corresponds to adding an interval with Hamiltonian
$\operatorname{diag}(1,0)$ at the left endpoint
\cite[Corollary~3.3]{KL}.
The condition $\int_0^x h_{22}(t)\,dt>0$ $(x>0)$
at the left endpoint in Theorem~\ref{thm_404} excludes such an interval.
Thus the first row of $U_a$ is also uniquely determined.
\end{proof}

The Hermite--Biehler function obtained from the second row of the fundamental matrix
at $x=L_a$ is therefore $\mathcal E_a=\sqrt\pi(A_a-iB_a)$.
By Proposition~\ref{prop_502}, the de Branges space it generates
is precisely $\mathcal B_a$.

Finally, we prove that $L_a$ depends continuously on $a$.
Set
\[
 \ell_a:=\int_0^{L_a}h_{11}(x)\,dx,\qquad
 m_a:=\int_0^{L_a}h_{22}(x)\,dx,\qquad
 L_a=\ell_a+m_a.
\]
Later, $\ell_a$ and $m_a$ will be interpreted as the length and mass of a string.

\begin{lemma} \label{lem_506}
Assume that $\lambda_R>0$ and set $L_0=0$.
The map $a\mapsto L_a$ in Proposition~\ref{prop_504}
is a continuous, strictly increasing bijection from $[0,R]$ onto $[0,L_R]$,
and its inverse is continuous.
The following hold as identities of Stieltjes measures in the variable $a$:
\begin{align}
&da=\sqrt{\frac{d\ell_a}{dL_a}\frac{dm_a}{dL_a}}\,dL_a,
\qquad
a=\int_0^{L_a}\sqrt{\det H_{\rm tr}(x)}\,dx,
\label{eq_506}\\
&dB_a(z)=zA_a(z)\,dm_a,
\qquad
dA_a(z)=-zB_a(z)\,d\ell_a.
\label{eq_507}
\end{align}
Moreover, $(A_0,B_0)=(1,0)$.
\end{lemma}
\begin{proof}
Let $U(x,z)=(U_{jk}(x,z))_{j,k=1}^2$ be the common fundamental matrix
in Proposition~\ref{prop_504}.
By Proposition~\ref{prop_505},
$
 (U_{21}(L_a,z),U_{22}(L_a,z))
 =(-B_a(z),A_a(z))
$.
Set $\alpha(x)=iU_{21}(x,i)$ and $\beta(x)=U_{22}(x,i)$.
Since the Hamiltonian is diagonal, \eqref{eq_402} gives
\[
 \alpha'=h_{22}\beta,\qquad
 \beta'=h_{11}\alpha,\qquad
 (\alpha(0),\beta(0))=(0,1),
\]
where the derivatives are with respect to $x$.
These differential equations and $h_{11},h_{22}\ge0$ imply
$\alpha\ge0$, $\beta\ge1$, and
\[
 \alpha(x)\ge\int_0^x h_{22}(y)\,dy>0
 \qquad(x>0).
\]
Both $\alpha$ and $\beta$ are nondecreasing.
Propositions~\ref{prop_502} and \ref{prop_505} give
$\alpha(L_a)=\alpha_a$ and $\beta(L_a)=\beta_a$. 
Let $\Pi_a$ denote the orthogonal projection
from $\mathcal B_R$ onto $\mathcal B_a$.
By the reproducing property, 
$J_a(w,\cdot)=\Pi_a J_R(w,\cdot)$. 
Thus the strong continuity of $\Pi_a$
in Proposition~\ref{prop_501} implies that
$J_a(i,0)$, $J_a(i,i)$, and $J_a(0,0)$
depend continuously on $a$. 
Since $\alpha_a>0$, the sum $\alpha_a+\beta_a$ is continuous for $a>0$.

Let $0<\varepsilon<R$. For $L_\varepsilon\le x\le L_R$, we have
\[
 \begin{aligned}
 (\alpha+\beta)'(x)
 &=h_{22}(x)\beta(x)+h_{11}(x)\alpha(x)\\
 &\ge\min(1,\alpha(L_\varepsilon))
       (h_{11}(x)+h_{22}(x))\\
 &=\min(1,\alpha(L_\varepsilon))>0
 \end{aligned}
\]
a.e.
Hence $\alpha+\beta$ is continuous and strictly increasing on this interval,
and its inverse is Lipschitz continuous.
It follows that
\[
 L_a=(\alpha+\beta)^{-1}(\alpha_a+\beta_a)
\]
is continuous on $[\varepsilon,R]$.

On the other hand, Proposition~\ref{prop_505} and the integral equation
for the fundamental matrix give
\[
 \int_0^{L_a}h_{22}(x)\,dx
 =B_a'(0)=J_a(0,0)~\to~0
 \quad(a\downarrow0).
\]
The last limit follows from the strong convergence of the orthogonal projections
to the projection onto the zero space in Proposition~\ref{prop_501}.
Since $\int_0^x h_{22}(y)\,dy>0$ for every $x>0$,
we obtain $L_a\to0$.
Thus $a\mapsto L_a$ is continuous on $[0,R]$.
By Proposition~\ref{prop_504}, the map $a\mapsto L_a$ is strictly increasing,
so it is a bijection from $[0,R]$ onto $[0,L_R]$,
with continuous inverse.

We next prove the identities of measures.
Since $a\mapsto L_a$ is continuous and strictly increasing,
a change of variables gives 
\[
 \ell_b-\ell_a
 =\int_{L_a}^{L_b}h_{11}(x)\,dx
 =\int_{(a,b]}h_{11}(L_t)\,dL_t
\]
for $0\le a<b\le R$. 
The same argument applies to $m_a$, so we obtain
\[
 d\ell_a=h_{11}(L_a)\,dL_a,\qquad
 dm_a=h_{22}(L_a)\,dL_a
\]
as identities of Stieltjes measures.
By Theorem~\ref{thm_404},
\[
 a=\int_0^{L_a}
 \sqrt{h_{11}(x)h_{22}(x)}\,dx.
\]
Applying the same change of variables gives
\[
 da=\sqrt{h_{11}(L_a)h_{22}(L_a)}\,dL_a
 =\sqrt{\frac{d\ell_a}{dL_a}
             \frac{dm_a}{dL_a}}\,dL_a.
\]
This proves \eqref{eq_506}.

For a diagonal Hamiltonian, the integral equation associated with
\eqref{eq_402} also gives
\[
 U_{21}(x,z)
 =-z\int_0^x U_{22}(y,z)h_{22}(y)\,dy,
 \qquad
 U_{22}(x,z)
 =1+z\int_0^x U_{21}(y,z)h_{11}(y)\,dy.
\]
Restricting to $x=L_a$ and using
$
(U_{21}(L_a,z),U_{22}(L_a,z))
=(-B_a(z),A_a(z))$
yields \eqref{eq_507}.
Finally, $U(0,z)=I_2$ gives
$(A_0,B_0)=(1,0)$.
\end{proof}

\section{Explicit recovery of the Hamiltonian}\label{section_6}

In Section~\ref{section_5}, we obtained the generating functions of each $\mathcal B_a$
from the fundamental matrix at $L_a$ of a single canonical system
with a trace-normed Hamiltonian.
In this section, we establish the local Lipschitz continuity
of $\alpha_a$ and $\beta_a$ and use it to prove that
the Stieltjes measures $d\ell_a$ and $dm_a$
in Lemma~\ref{lem_506} are absolutely continuous
with respect to Lebesgue measure $da$.
We then recover the Hamiltonian explicitly,
using $a$ as the independent variable,
and prove Theorem~\ref{thm_2}.

\subsection{Explicit recovery of the Hamiltonian}

By the formulas in Proposition~\ref{prop_502} and
$J_a(0,i)=\overline{J_a(i,0)}=\alpha_a>0$, we have
\begin{equation}\label{eq_601}
 A_a(i)=\beta_a,\qquad B_a(i)=i\alpha_a .
\end{equation}
We prove the following result in the latter part of this section.

\begin{proposition} \label{prop_601}
If $\lambda_R>0$, then $\alpha_a$ and $\beta_a$ are Lipschitz continuous
on each $[r,R]$ ($0<r<R$).
\end{proposition}

We use the following equalities from Lemma~\ref{lem_506}:
\begin{align}
 d A_a(z)=-zB_a(z)\,d\ell_a,\qquad
 &d B_a(z)=zA_a(z)\,dm_a,\qquad
 (A_0,B_0)=(1,0),\label{eq_602}\\
&da=\sqrt{\frac{d\ell_a}{dL_a}
           \frac{dm_a}{dL_a}}\,dL_a.
 \label{eq_603}
\end{align}
Here \eqref{eq_603} expresses
$a=\int_0^{L_a}\sqrt{\det H_{\rm tr}(x)}\,dx$
as an equality between measures.
At this stage, we know only that $da\ll dL_a$.

\begin{theorem} \label{thm_602}
Assume that $R>0$ and $\lambda_R>0$, 
and let $U_a(z)=U(L_a,z)$ be the fundamental matrix in Proposition~\ref{prop_504}. 
Then, for each $z\in\C$, the map
$a\mapsto U_a(z)$ is absolutely continuous on $[0,R]$ and satisfies
\begin{equation} \label{eq_604}
 \partial_a U_a(z)J=zU_a(z)H(a)
\end{equation}
for $da$-a.e.\ $a$.
The Hamiltonian $H(a)$ is unique up to equality $da$-a.e.\ and is given by
\begin{equation} \label{eq_605}
 H(a)=\begin{pmatrix}\gamma(a)^{-1}&0\\0&\gamma(a)\end{pmatrix},
 \quad
 \gamma(a)=\frac{\alpha_a'}{\beta_a}>0,
 \quad
 \gamma(a)^{-1}=\frac{\beta_a'}{\alpha_a}. 
\end{equation}
Here differentiation is with respect to $a$. 
The functions $\gamma$ and $\gamma^{-1}$ are integrable on $[0,R]$
and essentially bounded on each $[r,R]$ ($0<r<R$).
Moreover,
\begin{equation} \label{eq_606}
 \partial_a A_a(z)=-z\gamma(a)^{-1}B_a(z),
 \qquad
 \partial_a B_a(z)=z\gamma(a)A_a(z)
\end{equation}
holds $da$-a.e.
In particular, $\det H(a)=1$ holds $da$-a.e.
\end{theorem}
\begin{proof}
By Propositions~\ref{prop_504} and \ref{prop_505},
we have $U_a(z)=U(L_a,z)$ and
$(U_{a,21},U_{a,22})=(-B_a,A_a)$.
Evaluating \eqref{eq_602} at $z=i$ and using
\eqref{eq_601}, we obtain
\begin{equation} \label{eq_607}
 d\alpha_a=\beta_a\,dm_a,\qquad
 d\beta_a=\alpha_a\,d\ell_a,\qquad
 \alpha_{0+}=0,\quad \beta_{0+}=1.
\end{equation}
In particular, $\alpha_a$ and $\beta_a$ are nondecreasing and
$\beta_a\ge1$.
For every $0<r<R$,
$\alpha_a$ and $\beta_a$ are bounded below by positive constants on $[r,R]$
and are Lipschitz continuous there by Proposition~\ref{prop_601}.
Thus
\begin{equation} \label{eq_608}
 d\ell_a=\frac{\beta_a'}{\alpha_a}\,da,\qquad
 dm_a=\frac{\alpha_a'}{\beta_a}\,da
\end{equation}
holds on $(0,R]$.
Both measures are finite on $[0,R]$ and assign zero mass to the origin, 
so \eqref{eq_608} also holds on $[0,R]$.
Writing $\ell'$ and $m'$ for the densities in \eqref{eq_608}, 
we have $dL_a=(\ell'(a)+m'(a))\,da$. 
Since \eqref{eq_603} gives $da\ll dL_a$, it follows that
$\ell'(a)+m'(a)>0$ for $da$-almost every $a$.
Consequently,
\[
\frac{d\ell_a}{dL_a}
=\frac{\ell'(a)}{\ell'(a)+m'(a)},
\qquad
\frac{dm_a}{dL_a}
=\frac{m'(a)}{\ell'(a)+m'(a)}
\]
hold $dL_a$-almost everywhere.
Substituting these identities into \eqref{eq_603}, we obtain
\[
da
=\frac{\sqrt{\ell'(a)m'(a)}}{\ell'(a)+m'(a)}\,dL_a
=\sqrt{\ell'(a)m'(a)}\,da.
\]
Hence $\ell'(a)m'(a)=1$ for $da$-almost every $a$. 
In particular, both densities are positive $da$-almost everywhere.
Now set
\[
 \gamma(a):=m'(a)=\frac{\alpha_a'}{\beta_a}
\]
to obtain
$
 \gamma(a)^{-1}=\ell'(a)=\beta_a'/\alpha_a 
$
and hence \eqref{eq_605}.
Furthermore,
\[
 \int_0^R\bigl(\gamma(a)^{-1}+\gamma(a)\bigr)\,da
 =\ell_R+m_R=L_R<\infty
\]
shows that $\gamma$ and $\gamma^{-1}$ are integrable on $[0,R]$.
They are essentially bounded on each $[r,R]$ by Proposition~\ref{prop_601}
and the positive lower bounds for $\alpha_a,\beta_a$.

Next, we make the change of variable $x=L_a$
in the integral equation for the fundamental matrix
of the canonical system with the trace-normed
Hamiltonian $H_{\rm tr}$ constructed in
Proposition~\ref{prop_504}.
By Lemma~\ref{lem_506}, the measures corresponding to
$h_{11}(x)\,dx$ and $h_{22}(x)\,dx$
are $d\ell_a$ and $dm_a$, respectively.
Thus
\[
\begin{aligned}
U_a(z)
&=I_2-z\int_0^{L_a}U(x,z)H_{\rm tr}(x)J\,dx \\
&=I_2-z\int_{(0,a]}U_b(z)
   \begin{pmatrix}d\ell_b&0\\0&dm_b\end{pmatrix}J \\
&=I_2-z\int_0^a U_b(z)H(b)J\,db,
\end{aligned}
\]
where the last equality follows from
\eqref{eq_605} and \eqref{eq_608}.
Since $H$ is integrable on $[0,R]$,
the map $a\mapsto U_a(z)$ is absolutely continuous
on $[0,R]$ for each $z\in\C$ and satisfies \eqref{eq_604}.
Substituting \eqref{eq_608} into \eqref{eq_507}
also gives \eqref{eq_606}.
Since $\gamma(a)^{-1}\gamma(a)=1$,
we also have $\det H(a)=1$ for $da$-almost every $a$.

Finally, we prove that the Hamiltonian is unique.
If another Hamiltonian $\widetilde H$ satisfies \eqref{eq_604}
for the same fundamental matrix $U_a(z)$, then
subtracting the two differential equations at $z=i$ gives 
$U_a(i)(H(a)-\widetilde H(a))=0$ 
for $da$-almost every $a$.
Since $\det U_a(i)=1$ by the general property
of fundamental matrices of canonical systems,
$U_a(i)$ is invertible.
Therefore $\widetilde H=H$ $da$-almost everywhere. 
\end{proof}

\subsection{Proof of Theorem~\ref{thm_2}}

Assume that $Q_W^a>0$ for every $0<a<R$.
For any $0<b<R$, we have $\lambda_b>0$,
so Theorem~\ref{thm_602} applies with $R=b$.
Since $JH(a){}^tJ
=\operatorname{diag}(\gamma(a),\gamma(a)^{-1})$,
equation~\eqref{eq_606} agrees with \eqref{eq_105}.
Substituting $\alpha_a=J_a(i,0)$ and
$\beta_a=J_a(i,i)/J_a(i,0)$ into \eqref{eq_605}
gives \eqref{eq_106}.
\hfill $\Box$

\subsection{Differences of reproducing kernels}

The Hamiltonian recovered above gives an integral expression for the difference between the reproducing kernels of $\mathcal B_{a_1}$ and $\mathcal B_{a_2}$.

\begin{corollary} \label{cor_603}
Let $0\le a_1<a_2\le R$ and $\lambda_R>0$.
When $a_1=0$, set $\mathcal B_0=\{0\}$ and take its reproducing kernel to be zero.
For every $w,z\in\C$, we have
\begin{align}
 J_{a_2}(w,z)-J_{a_1}(w,z)
 &=\int_{a_1}^{a_2}\!\left\{
   \gamma(b)\overline{A_b(w)}A_b(z)
    +\gamma(b)^{-1}\overline{B_b(w)}B_b(z)\right\}db .
 \label{eq_609}
\end{align}
This difference is the reproducing kernel of $\mathcal B_{a_2}\ominus\mathcal B_{a_1}$.
In particular, 
\begin{equation} \label{eq_610}
 J_a(0,0)=\int_0^a\gamma(b)\,db,\qquad
 \gamma(a)=\frac{d}{da}\langle \mathsf A_a^{-1}1,1\rangle_{L^2(-a,a)}\quad\text{a.e.}
\end{equation}
\end{corollary}
\begin{proof}
Differentiating the numerator of \eqref{eq_503} with respect to $a$
and substituting \eqref{eq_606}, we obtain
\[
\frac{\partial}{\partial a}
\bigl(\overline{A_a(w)}B_a(z)-\overline{B_a(w)}A_a(z)\bigr)
=(z-\bar w)\bigl(\gamma(a)\overline{A_a(w)}A_a(z)
+\gamma(a)^{-1}\overline{B_a(w)}B_a(z)\bigr)
\]
for $da$-a.e.\ $a$.
For $z\ne\bar w$, we therefore obtain \eqref{eq_609}
by using \eqref{eq_503} and integrating from $a_1$ to $a_2$.
The case $a_1=0$ follows in the same way from $(A_0,B_0)=(1,0)$ and $J_0=0$.
Since both sides are continuous in $w$ and $z$,
the identity also holds when $z=\bar w$.
The space $\mathcal B_{a_1}$ is a closed subspace of $\mathcal B_{a_2}$,
and the orthogonal projection onto it maps $J_{a_2}(w,\cdot)$
to $J_{a_1}(w,\cdot)$.
The difference between the two kernels is therefore
the reproducing kernel of $\mathcal B_{a_2}\ominus\mathcal B_{a_1}$.
Finally, setting $a_1=0$, $a_2=a$, and $w=z=0$ in \eqref{eq_609}
and using $A_b(0)=1$, $B_b(0)=0$, and \eqref{eq_202},
we obtain \eqref{eq_610}.
\end{proof}

\subsection{Preliminaries for the proof of Proposition~\ref{prop_601}}

The terms involving prime powers in the local Weil quadratic form
(\eqref{eq_102} and Lemma~\ref{lem_201})
can be expressed using translation operators.
We rescale the interval $(-a,a)$ to $(-1,1)$
and estimate the effect of varying the translation parameters with $a$.
This requires uniform bounds not only for the $L^2$ norms of functions such as $u_{\pm i,a}$,
but also for the total variation of their extensions by zero.
For the remainder of this section, fix $0<r<R$
and let $a\in[r,R]$.
Define $V_a:L^2(-a,a)\to L^2(-1,1)$ by
$(V_au)(t)=\sqrt a\,u(at)$. 
Set
\[
\mathcal T_a=V_a\mathsf A_aV_a^{-1}.
\]
Let $P$ denote multiplication by $\mathbf1_{(-1,1)}$, and set
$\tau_su(t)=u(t-s)$.
In what follows, functions on $(-1,1)$ are extended by zero,
$\|\cdot\|_p$ denotes the norm of $L^p(-1,1)$
$(p=1,2)$, and $\langle\cdot,\cdot\rangle_2$ denotes
the inner product of $L^2(-1,1)$.
We write $\operatorname{Var}$ for the total variation on $\R$. 

We first express $\mathcal T_a$ as the sum of
a nonnegative self-adjoint operator independent of $a$
and a bounded self-adjoint operator.
For $0<\eta<\varepsilon<1$, set
\[
 j_{\eta,\varepsilon}(s):=
\frac{\mathbf1_{\{\eta<|s|<\varepsilon\}}}{2|s|},
\qquad
S_{\eta,\varepsilon}u:=
\log(\varepsilon/\eta)u-P(j_{\eta,\varepsilon}*u).
\]
The quadratic forms associated with $S_{\eta,\varepsilon}$
increase monotonically as $\eta\downarrow0$.
Let $S_\varepsilon$ be the self-adjoint operator associated with their limit.
We also set
\[
 \mathcal P_a:=\sum_{2\le n\le e^{2R}}\frac{\Lambda(n)}{\sqrt n}
 P\bigl(\tau_{(\log n)/a}+\tau_{-(\log n)/a}\bigr),\qquad
 \mathcal K_a u:=P(k_a*u),
\]
where
\begin{align*}
 k_a(s)&=\frac{\mathbf1_{\{\varepsilon\le |s|\le2\}}}{2|s|}
       +\mathbf1_{\{|s|\le2\}}a\,r''(as),\\
 r''(s)&=-2\cosh(s/2)+\frac{e^{-s/2}}{1-e^{-2s}}-\frac1{2s}
 \quad(s>0),\qquad r''(-s)=r''(s),
\end{align*}
and $r''$ is the second derivative of the remainder term in
\cite[\S2.2, \S2.4]{Su26a}.

\begin{lemma} \label{lem_604}
$S_\varepsilon$ is nonnegative and self-adjoint, with form domain
$H^{\log}(-1,1)$.
For $r\le a\le R$, the identity
\begin{equation} \label{eq_611}
 \mathcal T_a=S_\varepsilon+c_a\mathrm{Id}-\mathcal P_a-\mathcal K_a,
 \qquad c_a=-\log(a\varepsilon)-\log(2\pi)-C_0
\end{equation}
holds on the common operator domain $\mathfrak D(S_\varepsilon)$.
Moreover, for all $a,b\in[r,R]$, we have
\begin{equation} \label{eq_612}
 \sup_{r\le a\le R}\operatorname{Var}k_a<\infty,\qquad
 \|k_a-k_b\|_{L^1(\R)}
 \ll_{r,R,\varepsilon}|a-b|.
\end{equation}
\end{lemma}
\begin{proof}
We first establish the properties of $S_\varepsilon$.
Extending $u\in L^2(-1,1)$ by zero, we have
\[
 \langle S_{\eta,\varepsilon}u,u\rangle_2
 =\frac14\int_{\eta<|s|<\varepsilon}
     \|\tau_su-u\|_{L^2(\R)}^2
     \frac{ds}{|s|}.
\]
Thus the quadratic forms increase monotonically as $\eta\downarrow0$.
Their limit is a nonnegative closed form,
and the associated operator $S_\varepsilon$ is nonnegative and self-adjoint.
By Plancherel's identity, the limit form equals
\[
 \frac1{2\pi}\int_{-\infty}^{\infty}
 \left(\int_0^\varepsilon
       \frac{1-\cos(ts)}s\,ds\right)
 |\widehat u(t)|^2\,dt.
\]
The inner integral is $\log|t|+O_\varepsilon(1)$
as $|t|\to\infty$, so the form domain is $H^{\log}(-1,1)$.
Moreover, if $u\in C_c^\infty(-1,1)$ is extended by zero, then
\[
 \|\tau_su-u\|_{L^2(\R)}
 \le |s|\|u'\|_{L^2(\R)},
\]
so $S_{\eta,\varepsilon}u$ converges in $L^2(-1,1)$
as $\eta\downarrow0$.
Writing $h$ for its limit and $q_\varepsilon$ for the limit form, 
we have 
$q_\varepsilon(u,f)=\langle h,f\rangle_2$ for $f\in H^{\log}(-1,1)$. 
By the definition of the operator associated with a closed form,
it follows that $u\in\mathfrak D(S_\varepsilon)$ and
$S_\varepsilon u=h$.

Next, we derive the formula for $\mathcal T_a$.
Let $u\in C_c^\infty(-1,1)$ and set $v=V_a^{-1}u$.
If we define
\[
 F(s):=(v*\widetilde v)(as)
      =\langle\tau_{-s}u,u\rangle_2,
\]
then $F(0)=\|u\|_2^2$ and $F$ is supported in $[-2,2]$.
The function in the integral in \eqref{eq_102} can be written as
\[
 \frac{e^{x/2}}{e^x-e^{-x}}
 =\frac{e^{-x/2}}{1-e^{-2x}}
 =\frac1{2x}+r''(x)+2\cosh(x/2)
 \qquad(x>0).
\]
Moreover, when the integral is truncated at the origin, we have
\[
 2\int_\delta^\infty
       \frac{e^{-x}}{1-e^{-2x}}\,dx
 =\log\coth\frac{\delta}{2}
 =-\log\delta+\log2+o(1)
 \qquad(\delta\downarrow0).
\]
Substituting these identities into \eqref{eq_102}
and setting $\delta=a\eta$, we obtain
\[
\begin{aligned}
 Q_W^a(v,v)
 &=\lim_{\eta\downarrow0}
 \left\{
 -\log(a\eta)\|u\|_2^2
 -\frac12\int_{\eta<|s|\le2}
       \frac{F(s)}{|s|}\,ds
 \right\}\\
 &\quad-\sum_{n\le e^{2a}}
       \frac{\Lambda(n)}{\sqrt n}
       \left\langle
       (\tau_{(\log n)/a}
        +\tau_{-(\log n)/a})u,u
       \right\rangle_2\\
 &\quad-\bigl(\log(2\pi)+C_0\bigr)
       \|u\|_2^2-a\int_{-2}^2 F(s)r''(as)\,ds .
\end{aligned}
\]
By the definition of $S_{\eta,\varepsilon}$,
\[
\begin{aligned}
&-\log(a\eta)\|u\|_2^2
-\frac12\int_{\eta<|s|\le2}\frac{F(s)}{|s|}\,ds\\
&\qquad=
\langle S_{\eta,\varepsilon}u,u\rangle_2
-\log(a\varepsilon)\|u\|_2^2
-\frac12\int_{\varepsilon\le|s|\le2}
\frac{F(s)}{|s|}\,ds .
\end{aligned}
\]
Letting $\eta\downarrow0$ and using the definition of $k_a$, we therefore obtain
\[
 Q_W^a(V_a^{-1}u,V_a^{-1}u)
 =\langle S_\varepsilon u,u\rangle_2
 +c_a\|u\|_2^2
 -\langle\mathcal P_a u,u\rangle_2
 -\langle\mathcal K_a u,u\rangle_2.
\]
The operators $\mathcal P_a$ and $\mathcal K_a$
are bounded and self-adjoint on $L^2(-1,1)$.
Since $C_c^\infty(-1,1)$ is a form core for
$H^{\log}(-1,1)$, as recalled in Section~\ref{section_2},
the preceding identity extends to the entire form domain.
Thus the closed quadratic form of $\mathcal T_a$
agrees with that of 
$S_\varepsilon+c_a\mathrm{Id}-\mathcal P_a-\mathcal K_a$.
By uniqueness of the self-adjoint operator associated
with a closed quadratic form, the two operators coincide.
Since the remaining terms are bounded self-adjoint operators,
they do not change the operator domain.
Consequently, \eqref{eq_611} holds on the common domain
$\mathfrak D(S_\varepsilon)$.

Finally, we estimate $k_a$.
The formula for $r''$ gives
\[
 r''(s)=-\frac74-\frac{s}{48}+O(s^2)
 \qquad (s\downarrow0).
\]
Thus, if we set $r''(0)=-7/4$, then
$r''$ is continuous on $\R$ and Lipschitz continuous on every bounded interval.
Let $M$ bound the absolute value of $r''$ on $[-2R,2R]$,
and let $K$ be a Lipschitz constant on that interval.
For $r\le a,b\le R$ and $|s|\le2$, we have
\[
 \begin{aligned}
 |a r''(as)-b r''(bs)|
 &\le |a-b|\,|r''(as)|
       +b\,|r''(as)-r''(bs)|\\
 &\le (M+2RK)|a-b|,
 \end{aligned}
\]
and hence
\[
 \|a r''(a\cdot)-b r''(b\cdot)\|_{L^1(-2,2)}
 \le4(M+2RK)|a-b|.
\]
The function $s\mapsto a r''(as)$ is also Lipschitz continuous
on $[-2,2]$, with Lipschitz constant at most $R^2K$.
Indeed, for $s,t\in[-2,2]$,
we have $|ar''(as)-ar''(at)|\le a^2K|s-t|
\le R^2K|s-t|$.
Including the jumps at the two endpoints after extension by zero, we obtain
\[
 \operatorname{Var}
 \bigl(\mathbf1_{\{|s|\le2\}}a r''(as)\bigr)
 \le4R^2K+2RM.
\]
The first term in $k_a$ is independent of $a$
and has finite total variation because $\varepsilon>0$.
This proves \eqref{eq_612}.
\end{proof}

In the next lemma, $BV_0(-1,1)$ denotes the space of functions on $(-1,1)$
whose extensions by zero belong to $BV(\R)$,
with norm given by the total variation of the extension.
We do not require their traces at the endpoints to vanish.

\begin{lemma} \label{lem_605}
Assume that $\lambda_R>0$.
For $f_a=V_a1$ or $f_a=V_ae_{-i}$, we have
\[
 \sup_{r\le a\le R}\operatorname{Var}(\mathcal T_a^{-1}f_a)<\infty.
\]
\end{lemma}
\begin{proof}
For an integrable function $u$ of bounded variation, we have
$\operatorname{Var}(Pu)\le\operatorname{Var}(u)$.
Indeed, the size of the jump at either endpoint is bounded by the variation
from that endpoint towards infinity outside the interval, where $u$ tends to $0$.
The operator norms of $\mathcal P_a$ on $BV_0(-1,1)$ and $L^2(-1,1)$ are therefore both bounded by
\[
 M_R=2\sum_{2\le n\le e^{2R}}\frac{\Lambda(n)}{\sqrt n}.
\]
Moreover, since
$\operatorname{Var}(k_a*u)\le\operatorname{Var}(k_a)\|u\|_1$,
the operators $\mathcal K_a:L^2(-1,1)\to BV_0(-1,1)$ are uniformly bounded.
For any $c>0$, the operator norm of $P(j_{\eta,\varepsilon}*\cdot)$
on $BV_0(-1,1)$ is at most $\log(\varepsilon/\eta)$.
The Neumann series therefore gives
$\|(S_{\eta,\varepsilon}+c\mathrm{Id})^{-1}\|_{BV_0\to BV_0}\le1/c$. 
By the monotone convergence theorem for closed quadratic forms,
$(S_{\eta,\varepsilon}+c\mathrm{Id})^{-1}f$
converges to $(S_\varepsilon+c\mathrm{Id})^{-1}f$
in $L^2(-1,1)$ as $\eta\downarrow0$.
Since the interval is bounded, the convergence also holds
in $L^1(-1,1)$.
Thus, for $f\in BV_0(-1,1)$,
the lower semicontinuity of total variation gives
\[
\operatorname{Var}\bigl((S_\varepsilon+c\mathrm{Id})^{-1}f\bigr)
\le
\liminf_{\eta\downarrow0}
\operatorname{Var}\bigl(
(S_{\eta,\varepsilon}+c\mathrm{Id})^{-1}f
\bigr)\le c^{-1}\operatorname{Var}(f).
\]
Hence $(S_\varepsilon+c\mathrm{Id})^{-1}$
is bounded on $BV_0(-1,1)$ with norm at most $1/c$. 
The operator norm on $L^2(-1,1)$ is also at most $1/c$,
since $S_\varepsilon\ge0$.

Choose $\varepsilon$ sufficiently small that $c_*:=\min_{r\le a\le R}c_a>M_R$.
Since $\lambda_a\ge\lambda_R$, the function $w_a=\mathcal T_a^{-1}f_a$ satisfies
$\|w_a\|_2\le\|f_a\|_2/\lambda_R$.
Rearranging \eqref{eq_611} gives
\[
 (S_\varepsilon+c_a\mathrm{Id}-\mathcal P_a)w_a=f_a+\mathcal K_aw_a .
\]
On both spaces, the norm of $(S_\varepsilon+c_a\mathrm{Id})^{-1}\mathcal P_a$
is bounded by $M_R/c_*<1$.
Thus the operator on the left is invertible by a Neumann series,
and its inverse has norm at most $(c_*-M_R)^{-1}$ on $BV_0(-1,1)$.
Under the embedding $BV_0(-1,1)\hookrightarrow L^2(-1,1)$,
this inverse agrees with the inverse on $L^2$.
The function $f_a$ is either $\sqrt a$ or $\sqrt a e^{-at}$,
and its extension by zero has uniformly bounded total variation.
The $L^2$ estimate also bounds the total variation of $\mathcal K_aw_a$ uniformly,
which proves the lemma.
\end{proof}

In the operator decomposition, only the bounded remainder terms and finitely many translations vary.
Lemma~\ref{lem_605} also allows us to estimate the effect of varying the translation parameters
when the operators act on $\mathcal T_a^{-1}f_a$.

\subsection{Proof of Proposition~\ref{prop_601}}

Let each of $f_a,g_a$ be either $V_a1$ or $V_ae_{-i}$.
These families depend smoothly on $a$ as $L^2(-1,1)$-valued functions,
and $\|\mathcal T_a^{-1}\|\le1/\lambda_R$.
Hence
\begin{align*}
\langle\mathcal T_a^{-1}f_a,g_a\rangle_2
 -\langle\mathcal T_b^{-1}f_b,g_b\rangle_2
&=\langle\mathcal T_a^{-1}f_a,g_a-g_b\rangle_2
+\langle(\mathcal T_a^{-1}-\mathcal T_b^{-1})f_a,g_b\rangle_2\\
&\quad
 +\langle\mathcal T_b^{-1}(f_a-f_b),g_b\rangle_2\\
&=\langle(\mathcal T_a^{-1}-\mathcal T_b^{-1})f_a,g_b\rangle_2
 +O_{r,R}(|a-b|)\\
&=\big\langle(\mathcal T_b-\mathcal T_a)\mathcal T_a^{-1}f_a,
                    \mathcal T_b^{-1}g_b\big\rangle_2
 +O_{r,R}(|a-b|). 
\end{align*}
Here, by Lemma~\ref{lem_604}, the operators
$\mathcal T_a$ and $\mathcal T_b$ have the same operator domain.
Thus the resolvent identity
\[
\mathcal T_a^{-1}-\mathcal T_b^{-1}
=\mathcal T_b^{-1}(\mathcal T_b-\mathcal T_a)
 \mathcal T_a^{-1}
\]
holds on $L^2(-1,1)$.
The last equality above also uses the self-adjointness
of $\mathcal T_b^{-1}$.

Put $w_a=\mathcal T_a^{-1}f_a$ and
$v_b=\mathcal T_b^{-1}g_b$.
By \eqref{eq_611}, we have
\[
\mathcal T_b-\mathcal T_a
=(c_b-c_a)\mathrm{Id}
-(\mathcal P_b-\mathcal P_a)
-(\mathcal K_b-\mathcal K_a).
\]
We estimate the contributions of these three terms
to $\langle(\mathcal T_b-\mathcal T_a)w_a,v_b\rangle_2$.
First, the uniform bound on the inverses gives
$\|w_a\|_2+\|v_b\|_2\ll_{r,R}1$. 
Since
\[
|c_b-c_a|+\|k_b-k_a\|_{L^1(\R)}
\ll_{r,R}|a-b|
\]
by \eqref{eq_611} and \eqref{eq_612},
Young's inequality yields
\[
\left|
\langle(c_b-c_a)w_a,v_b\rangle_2
\right|
+
\left|
\langle(\mathcal K_b-\mathcal K_a)w_a,v_b\rangle_2
\right|
\ll_{r,R}|a-b|.
\]
It remains to estimate the contribution of
$\mathcal P_b-\mathcal P_a$.
By its definition,
\[
(\mathcal P_b-\mathcal P_a)w_a
=\sum_{2\le n\le e^{2R}}
\frac{\Lambda(n)}{\sqrt n}
P\left(
\tau_{(\log n)/b}-\tau_{(\log n)/a}
+\tau_{-(\log n)/b}-\tau_{-(\log n)/a}
\right)w_a.
\]
For functions extended by zero outside $(-1,1)$,
the standard inequalities for functions of bounded variation give
\[
\|\tau_s w-\tau_t w\|_{L^1(\R)}
\le |s-t|\operatorname{Var}(w),
\qquad
\|w\|_{L^\infty(\R)}
\le\operatorname{Var}(w).
\]
Consequently,
\begin{align*}
\left|
\langle P(\tau_s-\tau_t)w_a,v_b\rangle_2
\right|
&\le
\|\tau_s w_a-\tau_t w_a\|_{L^1(\R)}
\|v_b\|_{L^\infty(\R)}\\
&\le
|s-t|\operatorname{Var}(w_a)
\operatorname{Var}(v_b).
\end{align*}
Lemma~\ref{lem_605} bounds the two variations
uniformly for $a,b\in[r,R]$.
Moreover, for each $n$ in the finite sum,
\[
\left|\frac{\log n}{a}-\frac{\log n}{b}\right|
\le\frac{\log n}{r^2}|a-b|.
\]
Applying the preceding estimate to both
the positive and negative translations, we obtain
\[
\left|
\langle(\mathcal P_b-\mathcal P_a)w_a,v_b\rangle_2
\right|
\ll_{r,R}|a-b|.
\]
This estimate also covers the case where
$a$ and $b$ lie on opposite sides of $(\log n)/2$.
Indeed, the sum defining $\mathcal P_a$ has the fixed range
$2\le n\le e^{2R}$, and
$P\tau_s w=0$ for $|s|\ge2$
when $w$ is supported in $[-1,1]$.
Thus the corresponding translation term vanishes
automatically before it enters the interval.
Combining these estimates, we conclude that
\[
\langle\mathcal T_a^{-1}f_a,g_a\rangle_2
-\langle\mathcal T_b^{-1}f_b,g_b\rangle_2
=O_{r,R}(|a-b|).
\]
By \eqref{eq_202} and since $V_a$ is unitary, we have
\[
 J_a(i,0)=\langle\mathcal T_a^{-1}V_ae_{-i},V_a1\rangle_2, 
\quad 
 J_a(i,i)=\langle\mathcal T_a^{-1}V_ae_{-i},V_ae_{-i}\rangle_2.
\]
The preceding estimates show that $\alpha_a$ and $J_a(i,i)$
are Lipschitz continuous on $[r,R]$.
Since $\alpha_a$ is bounded below by a positive constant on $[r,R]$,
the function $\beta_a=J_a(i,i)/\alpha_a$ is also Lipschitz continuous.
\hfill $\Box$

\section{Applications} \label{section_7}

\subsection{Local spectral measures and zeros of the zeta function}\label{section_7_1}

Let $0<a<\infty$ and $\lambda_a>0$.
In what follows, a local spectral measure means a spectral measure of
$g|_{(-2a,2a)}$.
By the proof of Proposition~\ref{prop_504}, these are precisely
the measures $\nu$ for which $F\mapsto F|_{\R}$ is an isometry
from $\mathcal B_a$ into $L^2(\nu)$,
and each $\mu_{a,\theta}$ is an orthogonal spectral measure.
Under RH, the spectral measure of the global function $g$ is
\[
 \mu_\xi=\sum_{\gamma\in\Gamma}m_\gamma\delta_\gamma
\]
\cite[Theorem~1.1, (1.7), (1.8)]{Su23}.
For every $0<a<\infty$, this is a nonorthogonal spectral measure of $g|_{(-2a,2a)}$.
Indeed, for every $b>a$, we have $\mathcal B_a\subsetneq\mathcal B_b$
by Proposition~\ref{prop_501},
and, under RH, Proposition~\ref{prop_202} shows that
$F\mapsto F|_{\R}$ is an isometry from $\mathcal B_b$ into $L^2(\mu_\xi)$.
Its restriction to $\mathcal B_a$ is the same map,
so the image of $\mathcal B_a$ is a proper closed subspace of $L^2(\mu_\xi)$.
Thus $\mu_\xi$ is a nonorthogonal spectral measure.

\begin{proposition} \label{prop_701}
For distinct phases $\theta,\theta'\in[0,2\pi)$,
the sets $\Gamma_{a,\theta}$ and $\Gamma_{a,\theta'}$ interlace:
between any two consecutive points of either set, there is exactly one point of the other.
For every local spectral measure $\nu$ and every $t\in\R$,
we have $\nu(\{t\})J_a(t,t)\le1$,
with equality only for the unique orthogonal spectral measure whose spectrum contains $t$.
If $r_1<r_2$ are consecutive points of $\Gamma_{a,\pi}$,
then $\nu((r_1,r_2))>0$ whenever $\nu\ne\mu_{a,\pi}$.
In particular, under RH, each open interval $(r_1,r_2)$ contains
at least one point of $\Gamma$, and
$m_\gamma J_a(\gamma,\gamma)<1$
for every $\gamma\in\Gamma$.
\end{proposition}
\begin{proof}
By \eqref{eq_203} and Proposition~\ref{prop_502},
\[
 W(a,\theta;z)=i\alpha_a(e^{i\theta}-1)A_a(z)
                  +\beta_a(1+e^{i\theta})B_a(z).
\]
Thus the spectrum for $0<\theta<2\pi$ is the zero set of $A_a-rB_a$,
where $r=(\beta_a/\alpha_a)\cot(\theta/2)$.
For $\theta=0$, it is the zero set of $B_a$.
The functions $A_a,B_a$ have no common real zero,
and a continuous choice of $\arg(A_a(t)+iB_a(t))$ has derivative
\[
 \frac{A_a(t)B_a'(t)-A_a'(t)B_a(t)}{A_a(t)^2+B_a(t)^2}
 =\frac{J_a(t,t)}{A_a(t)^2+B_a(t)^2}>0.
\]
Each spectrum consists of the points where this argument takes a prescribed value modulo $\pi$.
Hence two distinct spectra have no common point, and between any two consecutive points of either,
there is exactly one point of the other.
Each $t\in\R$ belongs to exactly one $\Gamma_{a,\theta}$.

Computing the norm of the reproducing kernel $J_a(t,\cdot)$ in $L^2(\nu)$ gives
\[
 J_a(t,t)=\int_{-\infty}^{\infty}|J_a(t,x)|^2\,d\nu(x)\ge J_a(t,t)^2\nu(\{t\}).
\]
If equality holds, then $\nu$ is concentrated on $t$ and the real zeros of $J_a(t,\cdot)$.
By \eqref{eq_503}, these points form the set $\Gamma_{a,\theta}$ containing $t$.
Computing the norm of $J_a(x,\cdot)$ at each point $x$ of this set,
and using the orthogonality of the kernels at distinct points, yields
$\nu(\{x\})=1/J_a(x,x)$.
Conversely, these equalities hold for $\mu_{a,\theta}$.

Finally, for consecutive zeros $r_1,r_2$ of $A_a$,
the kernels $J_a(r_1,\cdot)$ and $J_a(r_2,\cdot)$ are orthogonal,
and on the real axis we have
\[
 J_a(r_1,x)J_a(r_2,x)
 =\frac{B_a(r_1)B_a(r_2)A_a(x)^2}{(x-r_1)(x-r_2)}
 \qquad (x\ne r_1,r_2).
\]
This product vanishes at the endpoints $r_1,r_2$,
and outside the interval it has a fixed strict sign wherever $A_a(x)\ne0$.
If $\nu((r_1,r_2))=0$, then, since the integral of the product with respect to $\nu$ vanishes,
the measure $\nu$ must be concentrated on the zeros of $A_a$.
The same computation of the kernel norms as above gives $\nu=\mu_{a,\pi}$.
Under RH, applying these results to $\nu=\mu_\xi$
and using that this measure is nonorthogonal proves the two conclusions.
\end{proof}

\subsection{The extension parameter for the global screw function}

The function $g$ associated with the Riemann zeta function,
defined in Section~\ref{section_5_2}, is a screw function
on the whole real axis if and only if RH holds
\cite[Theorem~1.2]{Su23}.
Under RH, the Hamiltonian recovered from the reproducing kernel
determines the parameter that selects $g$ among the extensions
of the local screw function in Proposition~\ref{prop_503}.

\begin{proposition}\label{prop_702}
Assume RH and let $a>0$.
Let $\widetilde g$ be an extension of $g|_{(-2a,2a)}$,
and let $\Omega$ be its parameter in \eqref{eq_505}.
Then $\widetilde g$ agrees with $g$ on the whole real axis if and only if
\[
 \Omega(z)=\lim_{b\to\infty}
 \frac{(U_a(z)^{-1}U_b(z))_{12}}{(U_a(z)^{-1}U_b(z))_{22}}
 \quad(\Im z>0),
\]
where the numerator and denominator are the $(1,2)$ and $(2,2)$ entries 
of the matrix $U_a(z)^{-1}U_b(z)$, respectively. 
The limit is locally uniform in the upper half-plane.
It is the Weyl function of the canonical system obtained by restricting
the Hamiltonian in \eqref{eq_605} to $(a,\infty)$
and normalizing the fundamental matrix to $I_2$ at the left endpoint $a$.
\end{proposition}
\begin{proof}

Under RH, the function $g$ extends the local screw functions for all $a>0$,
and its Weyl function is $Q_\xi(z)=i(\xi'/\xi)(1/2-iz)$
\cite[(1.7), (1.8)]{Su23}.
By the proof of Proposition~\ref{prop_504},
the fundamental matrix of the canonical system on the half-line
corresponding to this Weyl function agrees with $U_b$ at each $b>0$.
By \eqref{eq_604}, the matrix $U_a(z)^{-1}U_b(z)$
is the value at the right endpoint $b$ of the fundamental matrix
for the canonical system on $[a,b]$,
normalized to $I_2$ at the left endpoint $a$.
Thus the quotient in the proposition is the Weyl function
of the canonical system with boundary parameter $0$ at the right endpoint $b$.
Since $\det H=1$,
we have $\int_a^b\operatorname{tr}H(c)\,dc\ge2(b-a)$.
The canonical system on $[a,\infty)$ is therefore in the limit-point case.
The change of variable
$x=\int_a^b\operatorname{tr}H(c)\,dc$
normalizes the Hamiltonian to $\operatorname{tr}H=1$.
By Proposition~\ref{prop_402}, the quotient in the proposition
converges at every $z\in\C_+$ to the Weyl function of this canonical system.
Since $H$ is positive definite $da$-a.e., the limit is finite on $\C_+$.
The Cayley transforms of these quotients are holomorphic with values in the unit disk. 
By Montel's theorem and the pointwise convergence established above, 
they converge locally uniformly. 
Applying the inverse Cayley transform gives the same conclusion for the original quotients.
By Proposition~\ref{prop_403},
substituting this Weyl function into \eqref{eq_505} gives $Q_\xi$.
Normalized extensions correspond bijectively to their Weyl functions,
and the associated linear fractional transformation is injective because $\det U_a=1$.
The condition in the proposition is therefore necessary and sufficient.
\end{proof}

\subsection{Spectral formulas for the length and mass of the string}\label{section_7_3}

A canonical system with a diagonal Hamiltonian corresponds
to the spectral problem for a Krein string
(\cite[\S8.4]{KL}, \cite[\S2.4]{Su25}).
A Krein string is described by its length $\ell>0$
and a nondecreasing right-continuous function $M(y)$
giving the cumulative mass up to position $y$ on the string.
With the tension normalized to $1$, its equation of motion is
\[
 du'(y)+s u(y)\,dM(y)=0.
\]
Here $s$ is the spectral variable, and $'$ denotes differentiation with respect to $y$.
In general, the mass distribution $M$ may include point masses.

In the correspondence between diagonal Hamiltonians and Krein strings,
we set $s=z^2$, where $z$ is the spectral variable of the canonical system.
Assume that $\lambda_a>0$ and apply this correspondence
to the diagonal Hamiltonian in Proposition~\ref{prop_504}.
Let $x$ be the coordinate of the canonical system normalized by $\operatorname{tr}H=1$,
and define the position $y$ on the string and its mass distribution $M$ by
\[
 y(x)=\int_0^x h_{11}(v)\,dv,\qquad
 M(y(x))=\int_0^x h_{22}(v)\,dv
 \qquad(0\le x\le L_a).
\]
The normalization $\operatorname{tr}H=1$ gives $y(x)+M(y(x))=x$.
For every $0\le x_1<x_2\le L_a$, we have
$\int_{x_1}^{x_2}\sqrt{\det H_{\rm tr}(x)}\,dx>0$,
so $y$ is strictly increasing.
Composing $x\mapsto\int_0^xh_{22}(v)\,dv$ with the continuous inverse of $y$
shows that $M$ is continuous, with $M(0)=0$. 
Thus there is no point mass at either endpoint.
The length and total mass of the string are
$y(L_a)=\ell_a$ and $M(\ell_a)=m_a$, respectively.
By Theorem~\ref{thm_602}, we have
$d\ell_b=\gamma(b)^{-1}\,db$ and
$dm_b=\gamma(b)\,db$ $(0<b\le a)$.
We now express this length and mass in terms of the reproducing kernel
obtained from the local Weil quadratic form
and the functions $u_{z,a}$ given by the Riesz representation theorem.

Let $\varphi(y,s),\psi(y,s)$ be the fundamental solutions of the string equation
determined by $\varphi(0,s)=1$, $\varphi'(0,s)=0$,
$\psi(0,s)=0$, $\psi'(0,s)=1$.
Here ${}^\prime$ denotes differentiation with respect to $y$.
By \cite[(8.9), (8.10)]{KL}, the fundamental matrix of the canonical system is
\[
 U(x,z)=
 \begin{pmatrix}
  \psi'(y,s)&z\psi(y,s)\\
  \varphi'(y,s)/z&\varphi(y,s)
 \end{pmatrix}
 \qquad(y=y(x),\ s=z^2).
\]
Since its second row at $x=L_a$ is $(-B_a(z),A_a(z))$,
the spectrum $\Gamma_{a,\pi}$ corresponds to a free left endpoint and a fixed right endpoint,
with boundary condition $\varphi(\ell_a,s)=0$,
whereas $\Gamma_{a,0}$ corresponds to two free endpoints,
with boundary condition $\varphi'(\ell_a,s)=0$.
In the latter case, the eigenfunction at $s=0$ is constant.

To express the length $\ell_a$ directly in terms of $u_{\pm i,a}$,
we define
\[
 f(t)=u_{i,a}(t)-u_{-i,a}(t)
 -\int_{-a}^t(u_{i,a}+u_{-i,a})(y)\,dy
 +2\alpha_a\mathbf1_{(0,a)}(t)
\]
for $-a<t<a$, 
and extend $f$ by zero to the real axis.

\begin{proposition}
If $\lambda_a>0$, then $f\in H^{\log}(-a,a)$ and
\begin{equation}\label{eq_701}
 \ell_a=
 \sum_{\gamma\in\Gamma_{a,\pi}}
 \frac{1}{\gamma^2J_a(\gamma,\gamma)}
 =\frac{Q_W^a(f,f)}{|W(a,\pi;0)|^2},
 \qquad
 m_a=J_a(0,0)
 =\langle\mathsf A_a^{-1}1,1\rangle_{L^2(-a,a)}.
\end{equation}
\end{proposition}
\begin{proof}
Let $s=z^2<0$.
For the Weyl function of the string with a fixed right endpoint,
set $\Omega=0$ in \eqref{eq_505}
and use the symmetric spectral measure of phase $\pi$
identified in the proof of Proposition~\ref{prop_505} to obtain
\[
 \frac{\psi(\ell_a,s)}{\varphi(\ell_a,s)}
 =\frac{U_{a,12}(z)}{zU_{a,22}(z)}
 =\int_{-\infty}^{\infty}
   \frac{d\mu_{a,\pi}(\lambda)}{\lambda^2-s}.
\]
Thus the spectral measure of the string is the pushforward of $\mu_{a,\pi}$
under $\lambda\mapsto\lambda^2$.
Since $\mu_{a,\pi}$ is symmetric and $\mu_{a,\pi}(\{0\})=0$,
the mass at $\gamma^2$ corresponding to 
$\gamma\in\Gamma_{a,\pi}$, $\gamma>0$, is $2/J_a(\gamma,\gamma)$.

At $s=0$, we have $\varphi(y,0)=1$ and $\psi(y,0)=y$,
so $\psi(\ell_a,s)/\varphi(\ell_a,s)\to\ell_a$
$(s\uparrow0)$.
Applying the monotone convergence theorem to the Stieltjes integral on the right gives
\[
 \ell_a=
 2\sum_{\substack{\gamma\in\Gamma_{a,\pi}\\\gamma>0}}
 \frac{1}{\gamma^2J_a(\gamma,\gamma)}
 =
 \sum_{\gamma\in\Gamma_{a,\pi}}
 \frac{1}{\gamma^2J_a(\gamma,\gamma)}.
\]
On the other hand, differentiating
$dB_b(z)=zA_b(z)\,dm_b$ in \eqref{eq_507}
at $z=0$ and using $A_b(0)=1$, $B_0(z)=0$
gives $m_a=\partial_z B_a(0)=J_a(0,0)$.
Moreover, \eqref{eq_202} yields
\[
 J_a(0,0)=\widehat{u_{0,a}}(0)
 =\langle\mathsf A_a^{-1}1,1\rangle_{L^2(-a,a)}.
\]
When both endpoints are free, we have
\[
 \frac{\psi'(\ell_a,s)}{\varphi'(\ell_a,s)}
 =\frac{U_{a,11}(z)}{zU_{a,21}(z)}
 =\int_{-\infty}^{\infty}
   \frac{d\mu_{a,0}(\lambda)}{\lambda^2-s}.
\]
Since $\varphi'(\ell_a,s)=-sm_a+O(s^2)$
and $\psi'(\ell_a,s)=1+O(s)$, 
the spectral measure assigns mass
$\mu_{a,0}(\{0\})=1/m_a$ to the origin.

Finally, we prove the formula involving $f$.
Reflection symmetry gives $\widehat{u_{\pm i,a}}(0)=\alpha_a$.
Thus the sum of the last two terms in the definition of $f$
vanishes at both endpoints and has distributional derivative
$-(u_{i,a}+u_{-i,a})+2\alpha_a\delta_0$.
This sum is a compactly supported function of bounded variation,
so its Fourier transform is $O((1+|z|)^{-1})$ on the real axis.
Since $u_{\pm i,a}\in H^{\log}(-a,a)$,
it follows that $f\in H^{\log}(-a,a)$.
Taking Fourier transforms and using \eqref{eq_203}, we obtain
\[
\begin{aligned}
z\widehat f(z)
=(z-i)\widehat{u_{i,a}}(z)
 -(z+i)\widehat{u_{-i,a}}(z)+2i\alpha_a
=W(a,\pi;z)-W(a,\pi;0).
\end{aligned}
\]
Hence $\widehat f(\gamma)=-W(a,\pi;0)/\gamma$
for $\gamma\in\Gamma_{a,\pi}$.
The local Parseval identity \eqref{eq_205}
and the series for $\ell_a$ established above
give the remaining equality in \eqref{eq_701}.
\end{proof}

\section{Asymptotic behavior of the fundamental matrix under RH}
\label{section_8}

Assume RH throughout this section. 
Then, by Proposition~\ref{prop_504} and Theorem~\ref{thm_602}, 
the canonical system extends over the entire
positive half-line in the variable $a$.
We determine the asymptotic behavior of its fundamental matrix and the
reproducing kernel as $a\to\infty$.
In the notation of Proposition~\ref{prop_203}, write
$\Gamma\cap(0,\infty)=\{\gamma_1<\gamma_2<\cdots\}$ and
$\Gamma_{a,\pi}\cap(0,\infty)=\{\gamma_{a,1}<\gamma_{a,2}<\cdots\}$.
By \cite[Theorem~1.1]{Su26b}, 
\begin{equation}\label{eq_801}
 \gamma_{a,n}\ge\gamma_n,\qquad
 \gamma_{a,n}\searrow\gamma_n\quad(a\to\infty).
\end{equation}
Define the real entire function
\[
 A_\infty(z):=\prod_{n=1}^\infty\left(1-\frac{z^2}{\gamma_n^2}\right).
\]
The product converges locally uniformly in $\C$ since
$\sum_{n\ge1}\gamma_n^{-2}<\infty$.
Each distinct positive zero is included only once. In particular,
$A_\infty(z)=\xi(1/2-iz)/\xi(1/2)$ if and only if all the zeros of $\xi$ are simple.
By \eqref{eq_502} and \cite[Theorem~1.2 and (5.1)]{Su26b}, we have
\begin{equation}\label{eq_802}
 A_a(z)=\prod_{n=1}^\infty\left(1-\frac{z^2}{\gamma_{a,n}^2}\right)
 ~\to~ A_\infty(z)
\end{equation}
locally uniformly in $\C$ as $a\to\infty$.
Set
\[
 q_a(z):=\frac{U_{a,12}(z)}{A_a(z)},\qquad
 Q_\xi(z):=i\frac{\xi'}{\xi}(1/2-iz).
\]

\begin{proposition}\label{prop_801}
Assume RH. The length and mass of the string in Section~\ref{section_7_3}
satisfy
\begin{equation}\label{eq_803}
 \ell_a\nearrow\ell_\infty
 :=\frac{\xi''(1/2)}{\xi(1/2)}
 =\sum_{\gamma\in\Gamma}\frac{m_\gamma}{\gamma^2}<\infty,
 \qquad
 m_a=J_a(0,0)\ge\frac{a^2}{\ell_\infty}~\to~\infty
\end{equation}
as $a\to\infty$.
In particular, the string corresponding to the system on the half-line
has finite length and infinite mass.
\end{proposition}
\begin{proof}
As shown in the proof of Proposition~\ref{prop_702}, $Q_\xi$ is the Weyl
function of the common canonical system on the half-line.
Since $L_a\ge 2a$, Proposition~\ref{prop_402}, with boundary parameter
$\Omega=0$, gives $q_a(z)\to Q_\xi(z)$ for $z\in\C_+$.
The convergence is locally uniform there, since the Cayley transforms
of the Nevanlinna functions $q_a$ are bounded by one and $Q_\xi$ is
finite on $\C_+$, by Montel's theorem as in the proof of
Proposition~\ref{prop_702}.

The integral representation used in the proof of the proposition in
Section~\ref{section_7_3} gives
\begin{equation}\label{eq_804}
 \frac{q_a(z)}{z}
 =\sum_{n=1}^\infty
 \frac{2}{J_a(\gamma_{a,n},\gamma_{a,n})(\gamma_{a,n}^2-z^2)}.
\end{equation}
Initially this identity holds for $z=iy$, $y>0$.
By \eqref{eq_701}, the series converges locally uniformly off its poles,
so the equality extends there by analytic continuation, with a removable
singularity at $z=0$ on the left.
The representation contains no additional constant in $q_a(z)/z$,
or equivalently no linear term in $q_a$, as also follows from the
representation of $q_{\widetilde g}$ in Section~\ref{section_4_4}.  
By \eqref{eq_801}, we have $\gamma_{a,n}\ge\gamma_1$ for all $n$,
and hence
\[
 \frac{1}{\gamma_{a,n}^2+y^2}
 \le\frac{1}{\gamma_{a,n}^2}
 \le\left(1+\frac{y^2}{\gamma_1^2}\right)
 \frac{1}{\gamma_{a,n}^2+y^2}
 \qquad(y>0).
\]
Thus, comparing the series in \eqref{eq_701} and \eqref{eq_804}
term by term, we obtain
\[
 \frac{q_a(iy)}{iy}\le\ell_a\le
 \left(1+\frac{y^2}{\gamma_1^2}\right)\frac{q_a(iy)}{iy}
 \qquad(y>0).
\]
By Theorem~\ref{thm_602}, $\ell_a=\int_0^a\gamma(b)^{-1}\,db$ is increasing.
Taking the limit as $a\to\infty$ with $y>0$ fixed,
we see from these inequalities that $\ell_\infty$ is finite.
Then, taking the limit as $y\downarrow0$,
we obtain $\ell_\infty=Q_\xi'(0)$, since $Q_\xi(0)=0$. 
The functional equation gives $\xi'(1/2)=0$, and therefore
$Q_\xi'(0)=\xi''(1/2)/\xi(1/2)$.
For the expression in terms of the zeros, the Hadamard factorization
of the even entire function $\xi(1/2-iz)$ of order one gives
\[
 \frac{\xi(1/2-iz)}{\xi(1/2)}
 =\prod_{n=1}^\infty\left(1-\frac{z^2}{\gamma_n^2}\right)^{m_{\gamma_n}}.
\]
The sum $\sum_n m_{\gamma_n}/\gamma_n^2$ converges by the usual
counting estimate for the zeros of $\xi$.
Therefore, logarithmic differentiation gives
\[
 Q_\xi(z)=2z\sum_{n=1}^\infty\frac{m_{\gamma_n}}{\gamma_n^2-z^2}
 \qquad(z\notin\Gamma),
\]
where the series converges locally uniformly away from the poles.
In particular,
$\ell_\infty=Q_\xi'(0)=2\sum_n m_{\gamma_n}/\gamma_n^2>0$.
Finally, since $m_a=\int_0^a\gamma(b)\,db$ by
Theorem~\ref{thm_602}, the Cauchy--Schwarz inequality gives
$a^2\le\ell_a m_a\le\ell_\infty m_a$.
Thus, we obtain \eqref{eq_803}.
\end{proof}

As a consequence of Proposition~\ref{prop_801}, the average of
the trace-normed Hamiltonian tends to a matrix of rank one.
Indeed, we have
\[
 \frac{1}{L_a}\int_0^{L_a}H_{\rm tr}(x)\,dx
 =\frac{1}{\ell_a+m_a}
 \begin{pmatrix}\ell_a&0\\0&m_a\end{pmatrix}
 ~\to~\begin{pmatrix}0&0\\0&1\end{pmatrix}
 \qquad(a\to\infty), 
\]
since $L_a=\ell_a+m_a$. 
We next describe the asymptotic behavior of the fundamental matrix
$U_a(z)$ and the reproducing kernel $J_a(w,z)$ as $a\to\infty$.

\begin{theorem}\label{thm_802}
Assume RH. With $A_\infty$ and $Q_\xi$ as above, we have
\begin{equation}\label{eq_805}
 U_a(z)\begin{pmatrix}m_a^{-1}&0\\0&1\end{pmatrix}
 ~\to~
 \begin{pmatrix}
  -zQ_\xi(z)A_\infty(z)&Q_\xi(z)A_\infty(z)\\
  -zA_\infty(z)&A_\infty(z)
 \end{pmatrix}
\end{equation}
locally uniformly in $\C$ as $a\to\infty$.
Here the products involving $Q_\xi$ are understood by removing the
singularities at the points of $\Gamma$.
The limiting matrix has rank one at every $z\in\C$.
Moreover, we have 
\begin{equation} \label{eq_806}
 \frac{J_a(w,z)}{J_a(0,0)}~\to~
 \overline{A_\infty(w)}A_\infty(z)
\end{equation}
locally uniformly for $(w,z)\in\C^2$. 
\end{theorem}
\begin{proof}
Write $F(z)=\xi(1/2-iz)$. Since $Q_\xi=-F'/F$,
its poles at the points $\gamma\in\Gamma$ are simple and have residues
$-m_\gamma$. The zeros of $A_\infty$ are precisely the points of $\Gamma$
and are simple. Hence $Q_\xi A_\infty$ extends to an entire function with
\begin{equation}\label{eq_807}
 (Q_\xi A_\infty)(\gamma)=-m_\gamma A_\infty'(\gamma)\ne0
 \qquad(\gamma\in\Gamma).
\end{equation}

We first show that
\begin{equation}\label{eq_808}
 U_{a,12}(z)~\to~Q_\xi(z)A_\infty(z)
\end{equation}
locally uniformly in $\C$.
Multiplying \eqref{eq_804} by $zA_a(z)$ and using the product in
\eqref{eq_802}, we obtain
\[
 U_{a,12}(z)=z\sum_{n=1}^\infty
 \frac{2}{\gamma_{a,n}^2J_a(\gamma_{a,n},\gamma_{a,n})}
 \prod_{k\ne n}\left(1-\frac{z^2}{\gamma_{a,k}^2}\right).
\]
For $|z|\le r$ ($r>0$), each product is bounded by
$\prod_{k\ge1}(1+r^2/\gamma_k^2)$, by \eqref{eq_801}.
Thus \eqref{eq_701} gives locally uniform absolute convergence of the
series, which is valid also at the zeros of $A_a$. Moreover,
Proposition~\ref{prop_801} yields the bound
\[
 |U_{a,12}(z)|\le r\ell_\infty
 \prod_{k=1}^\infty\left(1+\frac{r^2}{\gamma_k^2}\right)
 \qquad(|z|\le r),
\]
uniformly in $a>0$.
On $\C_+$, we have $U_{a,12}=q_a A_a\to Q_\xi A_\infty$ by
\eqref{eq_802} and the convergence of the Weyl functions established
in the proof of Proposition~\ref{prop_801}.
The locally uniform bound and Montel's theorem, together with the identity
theorem, give \eqref{eq_808} on the whole plane.

Since the second row of $U_a$ is $(-B_a,A_a)$, the matrix equation
\eqref{eq_604} and the initial condition $U_0=I_2$ give
\[
 B_a(z)=z\int_0^a A_b(z)\gamma(b)\,db, \qquad 
 U_{a,11}(z)=1-z\int_0^a U_{b,12}(z)\gamma(b)\,db. 
\]
We divide by $m_a=\int_0^a\gamma(b)\,db$ and use weighted averages.
For either $V_b=A_b$ or $V_b=U_{b,12}$, let $V_\infty$ be its limit in
\eqref{eq_802} or \eqref{eq_808}, respectively.
On a compact set $K\subset\C$, the integral of
$|V_b(z)-V_\infty(z)|\gamma(b)$ over a fixed initial interval $[0,b_0]$
is uniformly bounded in $z\in K$, by the local integrability of $H$
and the integral equation for $U_b$.
The contribution of this interval to the weighted average tends to zero
because $m_a\to\infty$. The remaining contribution is bounded by
$\sup_{b\ge b_0,\,z\in K}|V_b(z)-V_\infty(z)|$, which tends to zero as
$b_0\to\infty$. Consequently,
\[
 \frac{B_a(z)}{m_a}~\to~ zA_\infty(z),\qquad
 \frac{U_{a,11}(z)}{m_a}~\to~-zQ_\xi(z)A_\infty(z)
\]
locally uniformly in $\C$. Together with \eqref{eq_802} and
\eqref{eq_808}, these limits prove \eqref{eq_805}.
Its right-hand side is the product of the column
${}^t(Q_\xi A_\infty,A_\infty)$ and the row $(-z,1)$.
By \eqref{eq_807}, the two entries of the column have no common zero.
Thus the rank is one at every point.

Finally, substitute the limits of $A_a$ and $B_a/m_a$ into
\eqref{eq_503} to obtain \eqref{eq_806} for $z\ne\bar w$.
To include the diagonal $z=\bar w$, put $v=\bar w$ and
$C_a=B_a/m_a$. The real entire functions $A_a,C_a$ satisfy
\[
 \frac{J_a(w,z)}{m_a}
 =A_a(v)\frac{C_a(z)-C_a(v)}{z-v}
  -C_a(v)\frac{A_a(z)-A_a(v)}{z-v}.
\]
For any entire function $G$, the quotient $(G(z)-G(v))/(z-v)$
can be written as $\int_0^1 G'(v+s(z-v))\,ds$, including the
case $z=v$ by continuity. 
Locally uniform convergence of the functions implies locally uniform
convergence of their derivatives. The last formula therefore gives
uniform convergence on compact sets of pairs $(w,z)$, with limit
$\overline{A_\infty(w)}A_\infty(z)$, as asserted.
\end{proof}

By \eqref{eq_805}, the norm of the first column of $U_a(z)$
tends to infinity as $a\to\infty$ for every fixed $z\ne0$.
Moreover, \eqref{eq_806} gives
$J_a(z,z)\sim m_a|A_\infty(z)|^2\to\infty$
for fixed $z\notin\Gamma$, whereas $J_a(\gamma,\gamma)$
remains bounded for each $\gamma\in\Gamma$
by Proposition~\ref{prop_701}.
Finally, for the generator $\mathcal E_a=\sqrt\pi(A_a-iB_a)$
from Proposition~\ref{prop_502}, 
we obtain
\[
 \frac{\mathcal E_a(z)}{m_a}\to-i\sqrt\pi\,zA_\infty(z)
\]
locally uniformly in $\C$.

\bigskip

\noindent\textbf{Acknowledgment.}
This work was supported by JSPS KAKENHI Grant Number JP23K03050.

%
\bigskip

\noindent\begin{minipage}{\linewidth}
Department of Mathematics, School of Science, Institute of Science Tokyo\\
2-12-1 Ookayama, Meguro-ku, Tokyo 152-8551, Japan\\[2pt]
Email: {\tt msuzuki@math.sci.isct.ac.jp}
\end{minipage}

\end{document}